\documentstyle[amsmath,amssymb,12pt]{article}
\renewcommand{\title}[1]{\leftline{\Large\bf #1}\par\medskip}
\renewcommand{\author}[1]{\medskip{\large #1}\par\medskip}

\newcommand{\om}{\omega}

\newcommand{\rom}[1]{{\rm #1}}

\allowdisplaybreaks

\makeatletter\@addtoreset{equation}{section}\makeatother

\begin{document}

\setcounter{page}{1} \setcounter{section}{0} \thispagestyle{empty}

\newtheorem{definition}{Definition}[section]
\newtheorem{remark}{Remark}[section]
\newtheorem{proposition}{Proposition}[section]
\newtheorem{theorem}{Theorem}[section]
\newtheorem{corollary}{Corollary}[section]
\newtheorem{lemma}{Lemma}[section]
\newtheorem{example}{Example}[section]

\newcommand{\skl}{\overset{(k)}{\diamondsuit}}
\newcommand{\D}{{\cal D}}
\newcommand{\N}{{\Bbb N}}
\newcommand{\C}{{\Bbb C}}
\newcommand{\Z}{{\Bbb Z}}
\newcommand{\R}{{\Bbb R}}
\newcommand{\Rp}{{\R_+}}
\newcommand{\eps}{\varepsilon}

\newcommand{\AS}{\operatorname{AS}}

\newcommand{\fii}{\varphi}
\newcommand{\Ker}{\operatorname{Ker}}
\newcommand{\supp}{\operatorname{supp}}
\newcommand{\la}{\langle}
\newcommand{\ra}{\rangle}
\newcommand{\const}{\operatorname{const}}
\newcommand{\ddGamma}{\overset{{.}{.}}{\Gamma}_X}

\newcommand{\ho}{\widehat\otimes}
\newcommand{\ot}{\otimes}

\renewcommand{\emptyset}{\varnothing}
\renewcommand{\tilde}{\widetilde}
\newcommand{\Formsg}{{\cal F}\Omega^{n+1}}
\newcommand{\Func}{{\cal FC}}

\newcommand{\di}{\partial}
\renewcommand{\div}{\operatorname{div}}

\begin{center}\LARGE \bf
DE RHAM COHOMOLOGY\\[2mm]\LARGE\bf OF CONFIGURATION SPACES\\[2mm]
\LARGE\bf WITH POISSON MEASURE \end{center}

\author{SERGIO ALBEVERIO}

\noindent{\sl Institut f\"{u}r Angewandte Mathematik,
Universit\"{a}t Bonn, Wegelerstr.~6, D-53115 Bonn, Germany;
\newline SFB 256, Univ.~Bonn, Germany;\newline SFB 237, Bochum--D\"usseldorf--Essen, Germany;
 \newline CERFIM (Locarno);
Acc.\ Arch.\ (USI), Switzerland;
\newline BiBoS, Univ.\ Bielefeld, Germany}\vspace{2mm}

\author{ALEXEI DALETSKII}

\noindent{\sl Institut f\"{u}r Angewandte Mathematik,
Universit\"{a}t Bonn, Wegelerstr.~6, D-53115 Bonn, Germany;
\newline SFB 256, Univ.~Bonn, Germany; \newline Institute of
Mathematics, Kiev, Ukraine;
\newline BiBoS, Univ.\ Bielefeld, Germany}\vspace{2mm}

\author{EUGENE LYTVYNOV}

\noindent{\sl Institut f\"{u}r Angewandte Mathematik,
Universit\"{a}t Bonn, Wegelerstr.~6, D-53115 Bonn, Germany;
\newline SFB 256, Univ.~Bonn, Germany;
\newline BiBoS, Univ.\ Bielefeld, Germany}

\begin{abstract}
\noindent The space $\Gamma_X$ of all locally finite
configurations in a
 Riemannian
manifold $X$ of infinite volume is considered. The deRham complex
of square-integrable differential forms over  $\Gamma_X$, equipped
with the Poisson measure, and the corresponding  deRham cohomology
are studied. The latter is shown to be   unitarily isomorphic to a
certain Hilbert tensor algebra generated by the $L^2$-cohomology
of the underlying manifold $X$.
\end{abstract}

\noindent 2000 {\it AMS Mathematics Subject Classification}.
  14F40, 58A10, 58A14, 60G57.

\newpage
\tableofcontents

\section{Introduction}

Let $\Gamma _X$ denote the space of all locally finite
configurations in a complete, stochastically complete, connected,
oriented Riemannian manifold $X$ of infinite volume. In this
paper, we define and study the deRham complex of square-integrable
differential forms over the configuration space $\Gamma_X$
equipped with the Poisson measure.

The  growing interest in geometry and analysis on the
configuration spaces can be explained by  the fact that these
naturally appear in different problems of   statistical mechanics
and quantum physics. In \cite{AKR-1, AKR0, AKR1}, an approach to
the configuration  spaces as  infinite-dimensional manifolds was
initiated. This approach  was motivated by the theory of
representations of diffeomorphism groups, see \cite{GGPS, VGG, I}
(these references as well as \cite{AKR1, AKR3} also contain
discussion of relations with quantum physics).  We refer the
reader to \cite{AKR2, AKR3, Ro, Lipsche} and references therein
for further discussion of analysis on the configuration  spaces
and applications.

On the other hand, stochastic differential geometry of
infinite-dimensional manifolds, in particular, their (stochastic)
cohomologies and related questions (Hodge--deRham Laplacians and
harmonic forms,  Hodge decomposition), has been  a very active
topic of research in recent years. It turns out that many
important examples of infinite-dimensional nonflat spaces (loop
spaces, product manifolds, configuration spaces) are naturally
equipped with probability measures (Brownian bridge, Gibbs
measures, Poisson measures). The geometry of these measures
interplays in a nontrivial way with the differential geometry of
the underlying  spaces themselves, and plays therefore a
significant role in their study. Moreover, in many cases the
absence  of a proper smooth manifold structure makes it  more
natural to work with $L^2$-objects (such as functions, sections,
etc.)\ on these infinite-dimensional  spaces, rather than to
define analogs of the smooth ones.

Thus, the concept of an $L^2$-cohomology has an important meaning
in this framework. The study of $L^2$-cohomologies for
finite-dimensional manifolds, initiated in \cite{Ati}, was a
subject of many works (whose different aspects  are treated  in
e.g.\ \cite{Dod, BrLes, ELR}, see also the review papers
\cite{Pan, Mat}). In the infinite-dimensional case,
 loop spaces have been most studied   \cite{JL, LRo, EL, Le}%
, the last two papers containing also a review of the subject. The
deRham complex on infinite product manifolds with Gibbs measures
(which appear in connection with problems of classical statistical
mechanics) was constructed in \cite{ADK1, ADK2} (see also
\cite{LeBe} for the case of the infinite-dimensional torus). We
should also mention the papers \cite{Shi, Ar1, Ar2, ArM, AK},
where the case of a flat Hilbert state space is considered (the
$L^2$-cohomological structure turns out to be nontrivial even in
this case  due to the existence of interesting measures on such a
space).

In \cite{ADL1, ADL2}, the authors started studying differential
forms over the infinite-dimensional space $\Gamma_X$, with $X$ as
above, and the corresponding Laplacians (of Bochner and deRham
type).

The structure of the present paper is as follows. Section~2 has an
introductory character. We recall the definition of the space
$L_\pi ^2\Omega ^n$ of differential forms over $\Gamma _X$ that
are square integrable with respect to the Poisson measure $\pi$,
and the construction of the unitary isomorhism
\[
I^n:L_\pi ^2\Omega ^n\rightarrow L_\pi ^2(\Gamma _X)\otimes \left[ %
\bigoplus_{m=1}^n L^2\Psi _{\mathrm sym}^n(X^m)\right],
\]
given in \cite{ADL2}. Here, $L_\pi ^2(\Gamma _X)$ is the space of
square-integrable functions over
$\Gamma _X$ and \linebreak $L^2\Psi _{\mathrm sym}^n(X^m)$ is
 a space of square-integrable $%
n$-forms over $X^m$ which satisfy some additional conditions.

We consider only the case of the Poisson measure with intensity
given by the Riemannian volume of $X$, which, according to
\cite{AKR1},  can be thought of as the volume measure on the
configuration space $\Gamma_X$, in the sense that the natural
lifting of the gradient and divergence on the underlying manifold
$X$ become dual operators.

In Section 3, we define the $L^2$-deRham complex over $\Gamma _X$
and the corresponding spaces ${\cal H}_\pi^{(n)}$ of (reduced)
$L^2$-cohomologies. We introduce the Hodge--deRham Laplacian ${\bf
H}^{(n)}$ acting in $L_\pi ^2\Omega ^n$ and study the space ${\bf
K}^{(n)}:=\operatorname{Ker}{\bf H}^{(n)}$ of harmonic forms. We
show, in particular, that ${\bf H}^{(n)}$ can be expressed, under
the action of the isomprphim $I^n$, in terms of the Laplacian
operator on functions on $\Gamma _X$ and the Hodge--deRham
Laplacians $H^{(n,m)}$ acting respectively in the spaces $L^2\Psi
_{\mathrm sym}^n(X^m)$. The application of the fact \cite{AKR1}
that the Dirichlet form of the Poisson measure is irreducible
gives us the possibility to express the harmonic forms on $\Gamma
_X$ in terms of harmonic forms on $X$. Our main result here is the
construction of the isomorphism
\[
\bigoplus_{n=0}^\infty{\bf K}^{(n)}\simeq {\cal A}_{\mathrm
sym}({\cal K}^{(1)},\dots,{\cal K}^{(\operatorname{dim}X)}),
\]
where ${\cal A}_{\mathrm sym}({\cal K}^{(1)},\dots,{\cal
K}^{(\operatorname{dim}X)})$ is a
supercommutative Hilbert tensor algebra generated by the spaces ${\cal K}%
^{(m)}:=\operatorname{Ker}H^{(m)}$, $H^{(m)}$ denoting the
Hodge--deRham Laplacian in the $L^2$-space of $m$-forms on $X$,
$m=1,\dots,\operatorname{dim}X$. The spaces ${\bf K}^{(n)}$ appear
to be  finite-dimensional, provided so are all the ${\cal
K}^{(m)}$ spaces. Using the weak Hodge--deRham decomposition, we
identify the spaces of harmonic forms with the spaces ${\cal
H}_\pi ^{n}$ of (reduced) $L^2$-cohomologies. In the case where
$\beta_m{:=}\operatorname{dim} {\cal K}^{(m)}<\infty$,
$m=1,\dots,\operatorname{dim} X$, we give an explicit formula for
the dimension $b_n$ of ${\cal H}_\pi ^{(n)}$:
\begin{equation*}\label{rdtrd}
b_n= \sum_{m=1}^n\,\sum_{1\le
k_1<\dots<k_m\le\operatorname{dim}X}\sum_{\begin{gathered}\scriptstyle
s_1,\dots,s_m\in\N:\\ \scriptstyle s_1 k_1+\dots+ s_m
k_m=n\end{gathered} }\beta _{k_1}^{(s_1)} \dotsm\beta
_{k_m}^{(s_m)},
\end{equation*}
where
\begin{equation*}\label{rfrtf}
\beta _k^{(s)}:=\begin{cases}\binom{ \beta _k} { s },&
k=1,3,\dots,\\ \binom{ \beta _k+s-1}{s} ,& k=2,4,\dots\end{cases}
\end{equation*}
We remark that this formula has the following interesting
consequence: although the spaces ${\cal H}_\pi ^{(n)}$ can be, in
general, nontrivial for any $n\in {\Bbb N}$, they  vanish for $n$
big enough, provided the cohomologies of $X$ of the even order do.

Finally, let us outline some links and open problems related to
the subject of the present paper.

1. Homology and homotopy of the spaces of {\it finite\/}
configurations (as topological spaces) were studied by many
authors (see e.g.\ \cite{fad, BCT, GKY, Wu}). An intriguing
question is to understand the relation between the results of
these authors and our results.

2. Any differential form $W\in L_\pi ^2\Omega ^n$ defines an
antisymmetric $n$-linear $L^2_\pi(\Gamma_X)$-valued functional on
the Lie algebra $\operatorname{Vect}_0(X)$ of compactly supported
vector fields over $X$. On the other hand, there exists a natural
representation of $\operatorname{Vect}_0(X)$ in $L_\pi ^2(\Gamma
_X)$ generated by the action of the diffeomorphism group
$\operatorname{Diff}_0(X)$ on $\Gamma _X $ (see \cite{VGG, AKR1,
AKR3}). It seems that the $L^2$-cohomology of $\Gamma _X$ is
related to a cohomolgy of the Lie algebra
$\operatorname{Vect}_0(X)$ with coefficients in this
representation.

3. In the present paper, we consider the case of the Poisson
measure with intensity given by  Riemannian volume of $X$. This
approach can  easily be extended to the case of a more general
intensity measure (the corresponding Hodge--deRham Laplacian is
defined in \cite{ADL2}). An important problem is to consider the
case of a Gibbs  measure (for analysis and geometry on
configuration spaces equipped with Gibbs measures and their
relations to the statistical mechanics of continuous systems, see
\cite{AKR2} and the review paper \cite{Ro}). The corresponding
$L^2$-cohomologies could give, in this case, invariants of such
measures and related models of statistical mechanics.

A different approach to the construction of differential forms and
related objects over Poisson spaces, based on the ``transfer
principle'' from Wiener spaces, is proposed in \cite{Pr2}, see
also \cite{PPr} and \cite{Pr}.

It is a great pleasure to thank K.~D.~Elworthy, Y.~G.~Kondratiev,
P.~Malliavin,
 M.~R\"{o}ckner, and A.~Thalmaier for their interest in
this work and helpful discussions. The financial support of SFB
256,  DFG Research Projects AL 214/9-3 and 436 UKR 113/43, and
BMBF Project UKR-004-99 is gratefully acknowledged.

\section{Differential forms over a configuration space}

The aim of this section is to recall some definitions and known
facts concerning the differential   structure of a configuration
space and differential forms over it. For more details and proofs,
we refer the reader to \cite{AKR1, ADL1, ADL2}.

Let $X$ be a complete, stochastically complete, connected,
oriented, $C^\infty $  Riemannian manifold of infinite volume. Let
$d$  denote  the dimension of $X$, $\langle \cdot ,\cdot
\rangle _x$ the  inner product in the tangent space $T_xX$ to $%
X $ at a point $x\in X$. The associated norm will be denoted by
$|\cdot |_x $. Let  $\nabla ^X$ stand for the gradient on $X$.

The configuration space $\Gamma _X$ over $X$ is defined as the set
of all locally finite subsets (configurations) in $X$:
\begin{equation}
\Gamma _X:=\big\{ \,\gamma \subset X\mid |\gamma \cap \Lambda
|<\infty \text{ for each compact }\Lambda \subset X\,\big\} .
\nonumber
\end{equation}
Here, $|A|$ denotes the cardinality of a set $A$.

We can identify any $\gamma \in \Gamma _X$ with the positive,
integer-valued Radon measure
\begin{equation}
\sum_{x\in \gamma }\varepsilon _x\subset {\cal M}(X),  \nonumber
\end{equation}
where $\varepsilon _x$ is the Dirac measure with mass at $x$,
$\sum_{x\in \varnothing }\varepsilon _x:=$zero measure, and ${\cal
M}(X)$ denotes the set of all positive Radon measures on the Borel
$\sigma $-algebra ${\cal B} (X)$. The space $\Gamma _X$ is endowed
with the relative topology as a subset of the space ${\cal M}(X)$
with the vague topology, i.e., the weakest topology on $\Gamma _X$
with respect to which  all maps
\begin{equation}
\Gamma _X\ni \gamma \mapsto \langle f,\gamma \rangle
:=\int_Xf(x)\,\gamma (dx)\equiv \sum_{x\in \gamma }f(x)  \nonumber
\end{equation}
are continuous. Here, $f\in C_0(X)$($:=$the set of all continuous
functions on $X$ with compact support). Let ${\cal B}(\Gamma _X)$
denote the corresponding Borel $\sigma $-algebra.

The tangent space to $\Gamma _X$ at a point $\gamma $ is defined
as the Hilbert space
\begin{equation}
T_\gamma \Gamma _X{:=}L^2(X\to TX;d\gamma )=\bigoplus_{x\in \gamma
}T_xX.  \label{tg-sp1}
\end{equation}
The scalar product and the norm in $T_\gamma \Gamma _X$ will be denoted by $%
\langle \cdot ,\cdot \rangle _\gamma $ and $\left\| \cdot \right\|
_\gamma $, respectively. Thus, each $V(\gamma )\in T_\gamma \Gamma
_X$ has the form $V(\gamma )=(V(\gamma )_x)_{x\in \gamma }$, where
$V(\gamma )_x\in T_xX$, and
\begin{equation}
\| V(\gamma )\| _\gamma ^2=\sum_{x\in \gamma }|V(\gamma )_x|_x^2.
\nonumber
\end{equation}

Vector fields and first order differential forms on $\Gamma _X$
will be identified with sections of the bundle $T\Gamma _X$.
Higher order
differential forms will be identified with sections of the tensor bundles $%
\wedge ^n(T\Gamma _X)$ with fibers
\begin{equation}
\wedge ^n(T_\gamma \Gamma _X)%
=\wedge ^n\left( \bigoplus_{x\in \gamma }T_xX%
\right),   \label{tang-n}
\end{equation}
where $\wedge ^n({\cal H})$ (or ${\cal H}^{\wedge n}$) stands for
the $n$th antisymmetric tensor power of a Hilbert space ${\cal
H}$.
Thus, under a differential form $W$ of order $n$, $n\in {\Bbb N}$, over $%
\Gamma _X,$ we will understand a mapping
\begin{equation}
\Gamma _X\ni \gamma \mapsto W(\gamma )\in \wedge ^n(T_\gamma
\Gamma _X). \label{lkghf}
\end{equation}

We will now recall how to introduce a  covariant derivative of a
differential form $W\colon \Gamma_X\to \wedge^n(T\Gamma_X)$.

Let $\gamma \in \Gamma _X$ and $x\in \gamma $. By ${\cal
O}_{\gamma ,x}$ we will denote an arbitrary open neighborhood of
$x$ in $X$ such that ${\cal O}_{\gamma ,x}\cap(\gamma \setminus
\{x\})=\varnothing$. We  define  the mapping
\begin{equation}\notag
{\cal O}_{\gamma ,x}\ni y\mapsto W_x(\gamma ,y)%
\mbox{$:=$}%
W(\gamma _y)\in \wedge ^n(T_{\gamma _y}\Gamma _X),\qquad
\gamma_y:=(\gamma\setminus\{x\})\cup\{y\}. \label{sec1}
\end{equation}
This  is a section of the Hilbert bundle
\begin{equation}
\wedge ^n(T_{\gamma _y}\Gamma _X)\mapsto y\in {\cal O}_{\gamma
,x}. \label{bund1}
\end{equation}
The Levi--Civita connection on $TX$ generates in a natural way a
connection on this bundle. We denote by $\nabla _{\gamma ,x}^X$
the corresponding covariant derivative and use the notation
\begin{equation}
\nabla _x^XW(\gamma )%
\mbox{$:=$}%
\nabla _{\gamma ,x}^X\,W_x(\gamma ,x)\in T_xX\otimes \left( \wedge
^n(T_\gamma \Gamma _X)\right)  \nonumber
\end{equation}
if the section $W_x(\gamma ,\cdot )$ is differentiable at $x$.

We say that the form $W$ is differentiable at a point $\gamma $ if
for each $x\in \gamma $ the section $W_x(\gamma ,\cdot )$ is
differentiable at $x$, and
\[
\nabla ^\Gamma W(\gamma )%
\mbox{$:=$}%
(\nabla _x^XW(\gamma ))_{x\in \gamma }\in T_\gamma \Gamma
_X\otimes \left( \wedge ^n(T_\gamma \Gamma _X)\right) .
\]
The mapping {\
\begin{equation}\notag
\Gamma _X\ni \gamma \mapsto \nabla ^\Gamma W(\gamma )%
\mbox{$:=$}%
(\nabla _x^XW(\gamma ))_{x\in \gamma }\in T_\gamma \Gamma
_X\otimes \left( \wedge ^n(T_\gamma \Gamma _X)\right)
\end{equation}
will be called the covariant gradient of the form $W$.

Analogously, one can introduce higher order derivatives of a
differential form $W$, the $m$th derivative
$(\nabla^\Gamma)^{(m)}W(\gamma)\in(T_\gamma\Gamma_X)^{\otimes
m}\otimes (\wedge^n(T_\gamma\Gamma_X))$.

Let us note that, for any $\eta \subset \gamma $, the space
$\wedge
^n(T_\eta \Gamma _X)$ can be identified in a natural way with a subspace of $%
\wedge ^n(T_\gamma \Gamma _X)$. In this sense, we will use the
expression  $W (\gamma )=W (\eta )$ without additional
explanations.

A form $W:\Gamma_X\to\wedge^n(T\Gamma_X)$ is called local  if
there exists a compact $\Lambda=\Lambda(W)$ in $X$ such that
$W(\gamma)=W(\gamma_\Lambda)$ for each $\gamma\in\Gamma_X$.

Let ${\cal F}\Omega^n$ denote the set of all local, infinitely
differentiable forms $W:\Gamma_X\to\wedge ^n(T\Gamma_X)$ which are
polynomially bounded, i.e., for each $W\in{\cal F}\Omega^n$ there
exist a function $\varphi\in C_0(X)$ and $k\in\N$ such that
\begin{equation}\label{weghtse}\|W(\gamma)\|^2_{
\wedge^n(T_\gamma\Gamma_X)}\le\langle\varphi^{\otimes
k},\gamma^{\otimes k}\rangle\qquad\text{for all
}\gamma\in\Gamma_X.\end{equation} Below, we will give an explicit
construction of a class of forms from ${\cal F}\Omega^n$.

Our next goal is to give a description of the space of $n$-forms
that are square-integrable with respect to the Poisson measure.

Let $dx$ denote the volume measure on $X$, and let  $\pi $ denote
the Poisson measure on $\Gamma _X$ with intensity $dx$. This
measure is characterized by its Laplace transform
\begin{equation}
\int_{\Gamma _X}e^{\langle f,\gamma \rangle }\,\pi (d\gamma
)=\exp\left[ \int_X(e^{f(x)}-1)\,dx\right],\qquad f\in C_0(X).
\nonumber
\end{equation}
If $F:\Gamma_X\to\R$ is integrable with respect to $\pi$ and
local, i.e., $F(\gamma )=F(\gamma _\Lambda )$ for some compact
$\Lambda\subset X$, then one has
\begin{equation}
\int_{\Gamma _X}F(\gamma )\,\pi (d\gamma )=e^{-\operatorname{vol}
(\Lambda
)}\sum_{n=0}^\infty \frac 1{n!}\int_{\Lambda ^n}F(\{x_1,\dots ,x_n\})\,dx_1%
\dotsm dx_n.  \label{3.1}
\end{equation}

We define on the set ${\cal F}\Omega ^n$ the $L^2$-scalar product
with respect to the Poisson measure:
\begin{equation}
(W_1,W_2)_{L_{\pi }^2\Omega ^n}%
\mbox{$:=$}%
\int_{\Gamma _X}\langle W_1(\gamma ),W_2(\gamma )\rangle _{\wedge
^n(T_\gamma \Gamma _X)}\,\pi (d\gamma ).  \label{4.1}
\end{equation}
The integral on the right hand side of \eqref{4.1} is finite,
since the Poisson measure has all moments finite. Moreover,
$(W,W)_{L_{\pi }^2\Omega ^n}>0$ if $W$ is not identically zero.
Hence, we can define a Hilbert space
$
L_{\pi }^2\Omega ^n%
=L^2(\Gamma _X\to \wedge ^n(T\Gamma _X);\pi )  
$ as the completion of ${\cal F}\Omega^n$ with respect to the norm
generated by the scalar product (\ref{4.1}).


We will now give an isomorphic description of the space $L_{\pi
}^2\Omega ^n$ via the space $L_{\pi }^2(\Gamma _X):=L^2(\Gamma
_X\to {\Bbb {R};\pi )}$ and some special spaces of
square-integrable forms on $X^m$, $m=1,\dots ,n$.

We first need  some preparations. For $x_1,\dots,x_n\in X$,  the
space $T_{x_1}X\wedge T_{x_2}X\wedge \dots \wedge T_{x_n}X$ will
be understood as a subspace of the Hilbert space
$\big(T_{y_1}X\oplus T_{y_2}X\oplus\dots\oplus
T_{y_k}X\big)^{\otimes n},$ where $\{y_1,\dots,y_k\}$ is the set
of the different $x_j$'s, $j=1,\dots,n$. We remark that
\begin{equation}\label{kjkjk}\big(T_{y_1}X\oplus
T_{y_2}X\oplus\dots\oplus T_{y_k}X\big)^{\otimes
n}\simeq\big(T_{y_{\nu(1)}}X\oplus
T_{y_{\nu(2)}}X\oplus\dots\oplus T_{y_{\nu(k)}}X\big)^{\otimes
n},\qquad \nu\in S_k\end{equation} (where $\simeq$ means
isomorphism), and moreover $T_{x_1}X\wedge
T_{x_2}X\wedge\dotsm\wedge T_{x_n}X$ and $T_{x_{\sigma(1)}}X\wedge
T_{x_{\sigma(2)}}X\wedge\dots\wedge T_{x_{\sigma(n)}}X$,
$\sigma\in S_n$, coincide   as  subspaces of the space
\eqref{kjkjk}.

Let
\[
\widetilde{X}^m:=\big\{\,(x_1,\dots,x_m)\in X^m\mid x_i\ne
x_j\text{ if }i\ne j\,\big\}.
\] Then, for $(x_1,\dots,x_m)\in\widetilde{X}^m$, we evidently have
\begin{equation}
\wedge ^n(T_{(x_1,\dots
,x_m)}X^m)=\bigoplus_{\begin{gathered}{\scriptstyle{
0\le k_1,\dots,k_m\le d}} \\ \scriptstyle k_1+\dots+k_m=n \end{gathered}%
}(T_{x_1}X)^{\wedge k_1}\wedge \dots \wedge (T_{x_m}X)^{\wedge
k_m}. \label{n-forms}
\end{equation}
For a form $\omega: X^m\to\wedge^n(TX^m)$ and
$(x_1,\dots,x_m)\in\widetilde{X}^m$, we denote by $\omega
(x_1,\dots ,x_m)_{k_1,\dots ,k_m}$ the corresponding component of
$\omega (x_1,\dots ,x_m)$ in the decomposition (\ref{n-forms}).

}We introduce the set $\Psi _{\mathrm sym}^n(X^m)$ of smooth forms
$\omega\colon X^m\to\wedge^n(TX^m) $ which have compact support
and satisfy the following assumptions on $\widetilde X^m$:

\begin{enumerate}
\item[(i)]  $\omega (x_1,\dots ,x_m)_{k_1,\dots ,k_m}=0$ if $k_j=0$ for some
$j\in \{1,\dots ,m\}$;

\item[(ii)]  $\omega $ is invariant with respect to the action of the group $S_m$:
\begin{equation}
\omega (x_1,\dots ,x_m)=\omega (x_{\sigma (1)},\dots ,x_{\sigma
(m)})\qquad \text{for each }\sigma \in S_m.  \label{symmetric}
\end{equation}

\end{enumerate}

For example, let $f\in C_0^\infty (X^2)$ be antisymmetric and let
$v\colon X\to TX$ be  a smooth, compactly supported vector field
on $X$. Then, the form $\omega\colon X^2\to \wedge^2(TX^2)$ given
by
\begin{align*}\omega(x_1,x_2):&= f(x_1,x_2)v(x_1)\wedge v(x_2)+
f(x_2,x_1) v(x_2)\wedge v(x_1)\\ &=2f(x_1,x_2)v(x_1)\wedge v(x_2)
\end{align*} belongs to $\Psi_{\mathrm sym} ^2(X^2)$.

Let us denote by $L ^2\Psi_{\mathrm sym} ^n(X^m)$ the Hilbert
space obtained as the completion of $\Psi _{\mathrm sym}^n(X^m)$
with respect to  the $L^2$-norm determined the  measure
$dx_1\dotsm dx_m$.

 We will use the notation
\begin{equation}
{\Bbb T}_{\{x_1,\dots ,x_m\}}^{(n)}X^m:=\bigoplus_{\begin{gathered}{%
\scriptstyle{ 1\le k_1,\dots,k_m\le d}} \\ \scriptstyle
k_1+\dots+k_m=n
\end{gathered}}(T_{x_1}X)^{\wedge k_1}\wedge \dots \wedge (T_{x_m}X)^{\wedge
k_m}.  \label{tang-n2}
\end{equation}
By virtue of \eqref{tang-n}, we have
\begin{equation}
\wedge ^n(T_\gamma \Gamma
_X)=\bigoplus_{m=1}^n\bigoplus_{\{x_1,\dots ,x_m\}\subset \gamma
}{\Bbb T}_{\{x_1,\dots ,x_m\}}^{(n)}X^m.  \label{tang-n1}
\end{equation}
For $W\in {\cal F}\Omega ^n$, we denote by $W_m(\gamma )\in
\bigoplus_{\{x_1,\dots ,x_m\}\subset \gamma }{\Bbb T}_{\{x_1,\dots
,x_m\}}^{(
n)}X^m$ the corresponding component of $W(\gamma )$ in the decomposition %
\eqref{tang-n1}. Thus, for $\{x_1,\dots ,x_m\}\subset \gamma $,
$W_m(\gamma ,x_1,\dots ,x_m)$ is equal to the projection of
$W(\gamma )\in \wedge ^n(T_\gamma \Gamma _X)$ onto the subspace
${\Bbb T}_{\{x_1,\dots ,x_m\}}^{(n)}X^m$.

\begin{proposition}\label{sq-int} \rom{\cite{ADL2}}
 Setting\rom, for $W\in
L^2_\pi\Omega^n$\rom,
\begin{equation}\label{asfdgfsdd}(I^nW)(\gamma,x_1,\dots,x_m):=(m!)^{-1/2}\,
W_m(\gamma\cup\{x_1,\dots,x_m\},x_1,\dots,x_m),\qquad m=1,\dots,n,
\end{equation}
one gets the unitary operator $$I^n:
L^2_\pi\Omega^n\to\bigoplus_{m=1}^n L^2_\pi(\Gamma_X)\otimes
L^2\Psi_{\mathrm sym}^n(X^m).$$
\end{proposition}

\begin{remark}
{\rm Actually,  formula \eqref{asfdgfsdd} makes sense only for
$(x_1,\dots,x_m)\in \widetilde{X}^m$.  However, since the set
$X^m\setminus \widetilde{X}^m$ is of zero $dx_1\dotsm dx_m$
measure, this does not lead to a contradiction.}
\end{remark}

\noindent{\it Sketch of the proof}. That $I^n$ is an isometric
operator from $L^2_\pi\Omega^n$ into $\bigoplus_{m=1}^n
L^2_\pi(\Gamma_X)\otimes L^2\Psi_{\mathrm sym}^n(X^m)$ follows
from the definition of $L^2\Psi^n_{\mathrm sym}(X^m)$,
\eqref{tang-n2}--\eqref{asfdgfsdd},  and the generalized Mecke
identity:
\begin{multline}\int_{\Gamma_X}\sum_{\{x_1,\dots,x_m\}\subset\gamma}
f(\gamma,x_1,\dots,x_m)\,\pi(d\gamma)\\=(m!)^{-1}\,\int_{\Gamma_X}\int_{X^m}
f(\gamma\cup\{x_1,\dots,x_m\},x_1,\dots,x_m)\,dx_1\dotsm
dx_m\,\pi(d\gamma),\label{awgf}\end{multline} where
$f:\Gamma_X\times X^m\to\R$ is a measurable function for which at
least one of the integrals in  \eqref{awgf} exists (this formula
can be proved by a repeated application of the Mecke identity, see
\cite{Gena}).

Let ${\cal FC}^\infty_{\mathrm b}({\cal D},\Gamma_X)$ denote the
set of smooth cylinder functions that is defined in Appendix~A.
For $F\in{\cal FC}^\infty_{\mathrm b}({\cal D},\Gamma_X)$ and
$\omega\in\Psi^n_{\mathrm sym}(X^m)$, $m\in\{1,\dots,n\}$, we
define a form $W$ by setting
\begin{equation}\label{kikimora}
W_k(\gamma,x_1,\dots,x_k):=\begin{cases}0,& k\ne m,\\ (m!)^{1/2}
F(\gamma\setminus\{x_1,\dots,x_m\})\om(x_1,\dots,x_m),&k=m.\end{cases}\end{equation}
As easily seen, $W$ is a local, infinitely differentiable $n$-form
over $\Gamma_X$ such that, for some $\varphi\in C_0(X)$,
$\varphi\ge0$, $$\|W(\gamma)\|^2_{\wedge^n(T_\gamma\Gamma_X)}\le
\langle \varphi^{\otimes n},\gamma^{\otimes n
}\rangle\qquad\text{for all }\gamma\in\Gamma_X,$$ and hence we
have the inclusion $W\in{\cal F }\Omega^n$. Moreover,
\begin{equation}\label{5634348}(I^nW)(\gamma,x_1,\dots,x_k)=\begin{cases}0,& k\ne
m,\\F(\gamma)\om(x_1,\dots,x_m),&k=m,\end{cases}\end{equation} for
each $\gamma\in\Gamma_X$ and each $(x_1,\dots,x_m)\in\tilde X^m$
such that $\{x_1,\dots,x_m\}\cap\gamma=\varnothing$. Since
$\gamma$ is a set of zero $dx$ measure and since the linear span
of $F\otimes\omega$ with $F$ and $\om$ as above, is dense in
$L^2_\pi(\Gamma_X)\otimes L^2\Psi_{\mathrm sym}^n(X^m)$, we obtain
the desired result.\quad$\blacksquare$\vspace{2mm}

In what follows, we will denote by ${\cal D}\Omega ^n$ the linear
span of the forms defined by (\ref{kikimora}) with $m=1,\dots,n$.
As we already noticed in the proof of Proposition \ref{sq-int},
${\cal D}\Omega ^n$ is a subset of
${\cal F} \Omega ^n$ and is dense in $L_{\pi }^2\Omega ^n$%
.

\section{De Rham complex over a configuration space}

\subsection{Exterior differentiation and $L^2$-cohomologies}

For $n\in\N$, let ${\cal E}\Omega^n$ denote the subset of ${\cal
F}\Omega^n$ consisting of all forms $W\in{\cal F}\Omega^n$ such
that all derivatives of $W$ are polynomially bounded, that is, for
each $k\in\N$ there exist $\varphi\in C_0(X)$, $\varphi\ge0$, and
$l\in\N$ (depending on $W$) such that
\begin{equation}\label{0989454}
\|(\nabla^\Gamma)^{(k)}W(\gamma)\|^2_{(T_\gamma\Gamma_X)^{\otimes
k}\otimes\wedge^n(T_\gamma\Gamma_X)} \le \langle \varphi^{\otimes
l},\gamma^{\otimes l}\rangle\qquad\text{for all
}\gamma\in\Gamma_X,\end{equation} and additionally, for each fixed
$\gamma\in\Gamma_X$ and $r\in\N$, the mapping $$
(X\setminus\gamma)^r\cap \widetilde X{}^r\ni(x_1,\dots,x_r)\mapsto
W(\gamma+\eps_{x_1}+\dots+\eps_{x_r})\in\wedge^n(T_\gamma\Gamma_X\oplus
T_{x_1}X\oplus\dots\oplus T_{x_r}X )$$ extends to a smooth,
compactly supported form $$X^r\ni(x_1,\dots,x_r) \mapsto
\omega(x_1,\dots,x_r)\in \wedge^n(T_\gamma\Gamma_X\oplus
T_{x_1}X\oplus\dots\oplus T_{x_r}X ).$$ (Notice that the locality
of a form, together with the above condition of extension, will
automatically imply the infinite differentiability of the form.)

As easily seen, ${\cal D}\Omega^n$ is a subset of ${\cal
E}\Omega^n$, and so we get the following  chain of inclusions
$${\cal D}\Omega^n\subset{\cal E}\Omega^n\subset{\cal
F}\Omega^n.$$

Absolutely analogously, we define the set ${\cal E}\Omega^0$ of
all local, smooth functions $F\colon \Gamma_X\to\R$ which,
together with all their derivatives, are polynomially bounded. We
have ${\cal FC}_{\mathrm b}^\infty({\cal D},\Gamma_X)\subset {\cal
F}\Omega^0$ (see Appendix~A).

 We define linear operators
\begin{equation}\label{awertudrz}
{\bf d}_n\colon {\cal E}\Omega ^n\to {\cal E}\Omega ^{n+1},\qquad n\in {\Bbb %
Z}_+,\end{equation} by
\begin{equation}\label{awgfse}
({\bf d}_nW)(\gamma ):=(n+1)^{1/2}\,\operatorname{AS}_{n+1}(\nabla
^\Gamma W(\gamma )),
\end{equation}
where \begin{equation}\label{uzui}\operatorname{AS}_{n+1}\colon
(T_\gamma \Gamma _X)^{\otimes (n+1)}\to \wedge ^{n+1}(T_\gamma
\Gamma _X)\end{equation} is the antisymmetrization operator. (We
notice that the polynomial boundedness of the form ${\bf d}_nW$
and its derivatives follows from the corresponding boundedness of
$\nabla^\Gamma W$ and the fact that the norm of the operator
\eqref{uzui} for each $\gamma\in\Gamma_X$ is equal to one).

Let us now consider ${\bf d}_n$ as an operator acting from the
space $L_\pi ^2\Omega ^n$ into $L_\pi ^2\Omega ^{n+1}$. We denote
by ${\bf d}_n^*$ the adjoint operator of ${\bf d}_n$.

\begin{proposition}\label{guzftdrt}
 ${\bf d}_n^{*}$ is a densely defined operator from
$L^2_\pi\Omega^{n+1}$ into $ L_\pi^2\Omega^n$ with domain
containing ${\cal E}\Omega ^{n+1}$\rom.
\end{proposition}

\noindent {\it Proof}. Let $\gamma \in \Gamma _X$ and $x\in \gamma $ be fixed. Let $%
C^\infty ({\cal O}_{\gamma ,x}\rightarrow \wedge ^n(T_\gamma
\Gamma _X))$ denote the space of all smooth sections of the
Hilbert bundle (\ref{bund1}). We define the operator
\[
d_{x,n}:C^\infty ({\cal O}_{\gamma ,x}\to \wedge ^n(T_\gamma
\Gamma _X))\to C^\infty ({\cal O}_{\gamma ,x}\rightarrow \wedge
^{n+1}(T_\gamma \Gamma _X))
\]
whose action, in local coordinates  on the manifold $X$, is given
as follows:
\[
d_{x,n}\,\phi (y)\,h_1\wedge \dots \wedge h_n=(n+1)^{1/2}\,\nabla
^X\phi (y)\wedge h_1\wedge \dots \wedge h_n,
\]
$\phi \in C^\infty ({\cal O}_{\gamma ,x}\to\R)$, $h_k\in
T_{x_k}X$, $x_k\in \gamma$, $k=1,\dots,n$. It easily follows from
the definition
of ${\bf d}_n$ and $%
\nabla ^\Gamma $ that, for $W\in{\cal F}\Omega^n$,
\begin{equation}
({\bf d}_nW)(\gamma )=\sum_{x\in \gamma }d_{x,n}W_x(\gamma ,x).
\end{equation}

Analogously, we define the operator
\[
\delta _{x,n}:C^\infty ({\cal O}_{\gamma ,x}\rightarrow \wedge
^{n+1}(T_\gamma \Gamma _X))\rightarrow C^\infty ({\cal O}_{\gamma
,x}\rightarrow \wedge ^n(T_\gamma \Gamma _X))
\]
setting
\begin{multline}
\delta _{x,n}\,\phi (y)\,h_1\wedge \dots \wedge
h_{n+1}:=\\=-(n+1)^{-1/2}\,\sum_{i=1}^{n+1}(-1)^{i-1}\varepsilon_{x,x_i}
\langle \nabla ^X\phi (y),h_i\rangle _xh_1\wedge \dots \wedge
\check{h}_i\wedge \dots \wedge h_{n+1},\label{awetz}
\end{multline}
where $\phi \in C^\infty ({\cal O}_{\gamma ,x}\to\R)$, $h_k\in
T_{x_k}X$, $x_k\in \gamma$, $k=1,\dots,n+1$,
\[
\varepsilon_{x,x_i}:=\begin{cases}1,& x=x_i,\\ 0,& x\ne
x_i,\end{cases}
\]
and $\check{h}_i$ denotes the absence of $h_i$. We now set for $W\in {\cal %
E}\Omega ^{n+1}$%
\begin{equation}\label{oipdrt}
{\pmb \delta }_nW(\gamma )=\sum_{x\in \gamma }\delta
_{x,n}W_x(\gamma ,x).
\end{equation}
By using \eqref{weghtse},  \eqref{awetz}, and \eqref{oipdrt}, we
conclude that
\begin{equation}\notag\label{rdzgzu}
{\pmb \delta }_n\colon {\cal E}\Omega ^{n+1}\to  {\cal E}\Omega^n.
\end{equation}
Moreover, from \eqref{3.1} and the definition of ${\bf d}_n$ and
${\pmb \delta}_n$, we derive, for arbitrary $V\in{\cal F}\Omega^n$
and $W\in{\cal F}\Omega^{n+1}$,
\begin{equation}\notag\label{lkawht}\int_{\Gamma_X}(({\bf
d}_nV)(\gamma),W(\gamma))_{\wedge^{n+1}(T_\gamma\Gamma_X)}\,\pi(d\gamma)
= \int_{\Gamma_X}(V(\gamma),({\pmb
\delta}_nW)(\gamma))_{\wedge^n(T_\gamma\Gamma_X)}\,\pi(d\gamma),\end{equation}
which proves the proposition. \quad $\blacksquare $

\begin{corollary}
The operator ${\bf d}_n:L_\pi^2\Omega ^n\rightarrow L_\pi^2\Omega
^{n+1}$\rom, $\operatorname{Dom}{\bf d}_n={\cal E}\Omega^n$\rom,
is closable\rom.
\end{corollary}

We denote by $\bar{\bf d}_n$ the closure of ${\bf d}_n$. The space $%
Z^n:=\operatorname{Ker}\bar{\bf d}_n$ is then a closed subspace of
$L_\pi ^2\Omega ^n$. Let  $B^n$ denote the closure in $L_\pi
^2\Omega ^n$ of the subspace $\operatorname{Im}{\bf d}_{n-1}$ (of
course, $B^n=$the closure of $\operatorname{Im}\bar{\bf
d}_{n-1}$).

We obviously have
$
{\bf d}_n{\bf d}_{n-1}=0 $, which implies $$ \operatorname{Im}{\bf
d}_{n-1}\subset\operatorname{Ker}{\bf d}_n\subset Z^n.$$ Hence
$B^n\subset Z^n$ and \begin{equation}\label{pserdo}
 \bar{\bf d}_n\bar{\bf
d}_{n-1}=0.\end{equation}

Thus, we have the infinite complex
\[
 \cdots\stackrel{{\bf d}_{n-1}}{\longrightarrow }{\cal E}\Omega ^n\stackrel{{\bf d}_n%
}{\longrightarrow }{\cal E}\Omega ^{n+1}\stackrel{{\bf d}_{n+1}}{\longrightarrow }%
\cdots\, ,
\]
 and the associated Hilbert complex
\begin{equation}
\cdots\stackrel{\bar{\bf d}_{n-1}}{\longrightarrow }L_\pi ^2\Omega ^n\stackrel{%
\bar{\bf d}_n}{\longrightarrow }L_\pi ^2\Omega ^{n+1}\stackrel{\bar{\bf d}%
_{n+1}}{\longrightarrow }\cdots\, .  \label{complex}
\end{equation}
Our next goal is to study the (reduced) $L^2$-cohomologies of
$\Gamma _X$, that is, the homologies of the complex
(\ref{complex}). We set in a standard way
\[
{\cal H}_\pi ^n=Z^n/B^n,\qquad n\in\N.
\]

Below, we will introduce the Hodge--deRham Laplacian  acting in
the space $L_\pi ^2\Omega ^n$, and identify  ${\cal H}_\pi ^n$
with
the space of harmonic forms. This will give us a possibility to express $%
{\cal H}_\pi ^n$ in terms of the cohomology spaces of the initial
manifold $X$.

\subsection{Hodge--deRham Laplacian of the Poisson measure}

For $n\in {\Bbb N}$, we define a bilinear form ${\cal E}_\pi
^{(n)}$ on $L^2_\pi\Omega^n$ by
\begin{multline}
{\cal E}_\pi ^{(n)}(W_1,W_2):=\int_{\Gamma _X}\big[ \langle {\bf d}%
_nW_1(\gamma ),{\bf d}_nW_2(\gamma )\rangle _{\wedge
^{n+1}(T_\gamma \Gamma _X)}  \label{lklk} \\
\text{{}}+\langle {\bf d}_{n-1}^{*}W_1(\gamma ),{\bf d}%
_{n-1}^{*}W_2(\gamma )\rangle _{\wedge ^{n-1}(T_\gamma \Gamma _X)}\big]%
\,\pi (d\gamma ),
\end{multline}
where $W_1,W_2\in \operatorname{Dom}{\cal E}_\pi ^{(n)}:={\cal
E}\Omega ^n$. The function under the sign of integral in
\eqref{lklk} is polynomially bounded, so that the integral exists.

\begin{theorem}
\label{th5.1} For any $W_1,W_2\in {\cal E}\Omega ^n${\rm , } we
have
\[
{\cal E}_\pi ^{(n)}(W_1,W_2)=\int_{\Gamma _X}\langle {\bf H}%
^{(n)}W_1(\gamma ),W_2(\gamma )\rangle _{\wedge ^n(T\Gamma
_X)}\,\pi (d\gamma ).
\]
Here{\rm , } ${\bf H}^{(n)}={\bf d}_{n-1}{\bf d}_{n-1}^{*}+{\bf d}_{n}^{*}{\bf %
d}_n$ is an operator in the space $L_\pi ^2\Omega ^n$ with domain
$\operatorname{Dom}{\bf H}%
^{(n)}:={\cal E}\Omega ^n${\rm . } It can be represented as
follows{\rm : }
\begin{equation}
{\bf H}^{(n)}W(\gamma )=\sum_{x\in \gamma }H_x^{(n)}W(\gamma
)=\langle H_{\bullet }^{(n)}\,W(\gamma ),\gamma \rangle ,\qquad
W\in {\cal E}\Omega ^n, \label{5.7}
\end{equation}
where
\begin{equation}
H_x^{(n)}=d_{x,n-1}\delta_{x,n-1}+\delta_{x,n}d_{x,n}. \label{5.8}
\end{equation}
\end{theorem}

\noindent{\it Proof}. The statement follows from
\eqref{awertudrz},  (the proof of) Proposition~\ref{guzftdrt}, and
the equality $d_{x,n-1}\delta_{y,n-1}+\delta_{y,n} d_{x,n}=0$
holding for all $x,y\in\Gamma$, $x\ne y$.\quad
\vspace{2mm}$\blacksquare$

From Theorem~\ref{th5.1} we conclude that the bilinear form ${\cal E}%
_\pi ^{(n)}$ is closable in the space $L_\pi ^2\Omega ^n$. The
generator of its
closure (being actually the Friedrichs extension of the operator ${\bf H}%
^{(n)}$, for which we preserve the same notation) will be called
the Hodge--deRham
Laplacian on $\Gamma _X$ (corresponding to the Poisson measure $\pi $). By (%
\ref{5.7}) and (\ref{5.8}), ${\bf H}^{(n)}$ is the lifting of the
Hodge--deRham Laplacian on $X$.


For linear operators $A$ and $B$ acting in Hilbert spaces ${\cal H}$ and $%
{\cal K}$, respectively, we introduce an operator $A\boxplus B$ in ${\cal H}%
\otimes {\cal K}$ by
\begin{equation}
A\boxplus B%
\mbox{$:=$}%
A\otimes {\bf 1}+{\bf 1}\otimes B,\qquad \operatorname{Dom}(A\boxplus B):=%
\operatorname{Dom}(A)\otimes _{\mathrm a}\operatorname{Dom}(B),
\nonumber
\end{equation}
where $\otimes _a$ stands for the algebraic tensor product. If the
operators
$A$ and $B$ are closable, then so is $%
A\boxplus B$, and we will preserve   the same notation for its
closure.

Next, for operators $A_1,\dots ,A_n$ acting in Hilbert spaces ${\cal H}%
_1,\dots ,{\cal H}_n$, respectively, let $\bigoplus_{i=1}^nA_i$
denote the operator in $\bigoplus_{i=1}^n{\cal H}_i$ given by
\[
\left(\bigoplus_{i=1}^nA_i\right)(f_1,\dots ,f_n)=(A_1f_1,\dots
,A_nf_n),\qquad f_i\in \operatorname{Dom}(A_i).
\]

\begin{theorem}
\label{thonsa2} {\rm 1)} On ${\cal D}\Omega ^n$ we have
\begin{equation}
{\bf H}^{(n)}=(I^n)^{-1}\left[{\bf H}^{(0)}\boxplus \left(
\bigoplus _{m=1}^n H_{\mathrm sym}^{(n,m)}\right)\right]I^n,
\label{mn}
\end{equation}
where ${\bf H}^{(0)}$ is the Laplacian in the space
$L^2_\pi(\Gamma_X)$ \rom(see Appendix~\rom{A),} and
$H_{\mathrm sym}^{(n,m)}$ is the restriction of the Hodge--deRham Laplacian $%
H^{(n,m)}$ acting in the space $L^2\Omega ^n(X^m):=L^2(X^m\to
\wedge^n(T X^m);\,dx_1\dotsm dx_m)$ to the subspace $L^2\Psi
_{\mathrm sym}^n(X^m)${\rm . }

{\rm 2)\ }  ${\cal D}\Omega ^n$ is a domain of essential selfadjointness of  $%
{\bf H}^{(n)}${\rm , } and the equality~\eqref{mn} holds for the
closed
operators ${\bf H}^{(n)}$ and ${\bf H}^{(0)}\boxplus \big( %
\bigoplus_{m=1}^nH_{\mathrm sym}^{(n,m)}\big)$ {\rm (}where the
latter operator is
closed from its domain of essential selfadjointness $I^n({\cal D}\Omega ^n)$%
{\rm )}{\rm . }
\end{theorem}

\noindent {\it Proof}. This theorem was proved in \cite{ADL2} in a
more general setting. Here, we present a simplified version of
this proof adapted to our special case of the volume measure on
$X$.

1) Let $W\in {\cal D}\Omega ^n$ be given by formula
\eqref{kikimora}. Then, using Theorem~\ref{th5.1} and Appendix~A,
we get
\begin{gather}
\left({\bf H}^{(n)}W\right)_k(\gamma )=0\qquad \text{for }k\ne m,
\nonumber
\\
\left({\bf H}^{(n)}W\right)_m(\gamma,\bar x^m )=\left(\sum_{x\in \gamma }H_x^{(n)}W%
\right)_m(\gamma,\bar x^m )  \nonumber \\ =\left(\sum_{x\in \gamma
\setminus \{\bar x^m \}} H_x^{(n)}W\right)_m(\gamma,\bar x^m
)+\left(\sum_{x\in \{\bar x^m \}}H_x^{(n)}W\right)_m(\gamma ,\bar
x^m ) \nonumber
\\
=(m!)^{1/2} \left[\left(\sum_{x\in \gamma \setminus \{\bar x^m \}}
H_x F\right)(\gamma \setminus \{%
\bar x^m \})\omega (\bar x^m )
+F(\gamma \setminus \{\bar x^m \})
\left(\sum_{x\in \{%
\bar x^m \}}H_x^{(n)}\omega \right)(\bar x^m )\right] \nonumber
\\ =(m!)^{1/2}\left[\left({\bf H}^{(0)}F\right)(\gamma \setminus \{\bar x^m \})\omega
(\bar x^m )+F(\gamma \setminus \{\bar x^m \})\left(H_{\mathrm
sym}^{(n,m)}\omega \right)(\bar x^m )\right], \label{chyzh}
\end{gather}
where $\bar x^m:=(x_1,\dots,x_m)$, $\{\bar
x^m\}:=\{x_1,\dots,x_m\}$, and $\{\bar x^m \}\subset\gamma$.
(Notice that the Hodge--deRham Laplacian in the space
$L^2\Omega^n(X^m)$ leaves the set $\Psi _{\mathrm sym} ^n(X^m)$
invariant.) Therefore,
\begin{equation}
(I^n{\bf H}^{(n)}W)(\gamma ,\bar{x}^k)=\begin{cases}0,& k\ne m,\\
({\bf H}^{(0)}F) (\gamma)\omega(\bar{x}^m)+F(\gamma)(H_{\mathrm
sym}^{(n,m)}\omega)(\bar{x}^m),&k=m.\end{cases} \label{pyzh}
\end{equation}
Hence, by virtue of \eqref{5634348}, we get \[\left(\left[{\bf
H}^{(0)}\boxplus \left( \bigoplus _{m=1}^n H_{\mathrm
sym}^{(n,m)}\right)\right]I^n\right)
(\gamma ,\bar{x}^k)=\left(I^n{\bf H}^{(n)}W\right)(\gamma ,%
\bar{x}^k),\qquad k=1,\dots ,n,
\]
which proves \eqref{mn}.


2) Let $\Omega^n(X^m)$ denote the space of all smooth forms $\omega \colon %
X^m\to \wedge ^n(TX^m)$ with compact support, and
let $L^2\Omega^n_{\mathrm sym}(X^m)$ denote the subspace of $%
L^2\Omega^n(X^m)$ consisting of all  forms invariant with respect
to the action of the symmetric group $S_m$, i.e., the
forms $\omega \in L^2\Omega^n(X^m)$ for which the equality %
\eqref{symmetric} holds for a.a.\ $(x_1,\dots ,x_m)\in
\widetilde{X}^m$. Evidently, the orthogonal projection $P_m^n$
onto this subspace is given by the formula
\begin{equation}
(P_m^n\omega )(x_1,\dots ,x_m)=\frac 1{m!}\sum_{\sigma \in
S_m}\omega (x_{\sigma (1)},\dots ,x_{\sigma (m)})  \label{kkk1}
\end{equation}
and
\begin{equation}
P_m^n\Omega ^n(X^m)=\Omega _{\mathrm sym}^n(X^m),  \label{kkk2}
\end{equation}
where $\Omega _{\mathrm sym}^n(X^m)$ denotes the subspace of
$\Omega^n(X^m)$ consisting of all $S_m$-invariant forms.

It is known that the Hodge--deRham Laplacian $H^{(n,m)}$ in
$L^2\Omega^n(X^m)$ is essentially self-adjoint on $\Omega ^n(X^m)$
(e.g.~\cite{E3}). Then, the nonnegative definiteness of
$H^{(n,m)}$ yields that the set $(H^{(n,m)}+{\bf 1})\Omega
^n(X^m)$ is dense in $L^2\Omega ^n(X^m)$, see e.g.\
\cite[Section~10.1]{RS}. Therefore, the set $P_m^n(H^{(n,m)}+{\bf
1})\Omega ^n(X^m)$ is dense in $L^2\Omega^n_{\mathrm sym}(X^m)$.
But upon \eqref{kkk1} and \eqref{kkk2},
\[
P_m^n(H^{(n,m)}+{\bf 1})\Omega
^n(X^m)=(H^{(n,m)}P_m^n+P_m^n)\Omega ^n(X^m)=(H^{(n,m)}+{\bf
1})\Omega _{\mathrm sym}^n(X^m),
\]
which implies that the restriction $H_{\mathrm sym}^{(n,m)}$ of
the operator $H^{(n,m)}$ to the subspace
$L^2\Omega^n_{\mathrm sym}(X^m)$ is essentially self-adjoint on $%
\Omega _{\mathrm sym}^n(X^m)$.

Because $H_{\mathrm sym}^{(n,m)}$ acts invariantly on the subspace
$L^2\Psi_{\mathrm sym} ^n(X^m)$ and its orthogonal complement in
$L^2\Omega_{\mathrm sym} ^n(X^m)$, we conclude that $H_{\mathrm
sym}^{(n,m)}$ considered as an operator in $L^2\Psi_{\mathrm sym}
^n(X^m)$ is essentially self-adjoint on $\Psi _{\mathrm
sym}^n(X^m)$. Consequently, the operator
$\bigoplus_{m=1}^nH_{\mathrm sym}^{(n,m)}$ is
essentially self-adjoint on the direct sum of the sets $\Psi _{\mathrm sym}^n(X^m)$, $%
m=1,\dots ,n$.

Finally, remarking that the operator ${\bf H}^{(0)}$ is
essentially self-adjoint on ${\cal FC}_{\mathrm b}^\infty({\cal D
},\Gamma_X)$ (\cite[Theorem~5.3]{AKR1}, see also Appendix~A), we
conclude from the theory of operators admitting separation of
variables (e.g.\ \cite[Ch.~6]{B})
that $I^n({\cal D}%
\Omega ^n)$ is a domain of essential self-adjointness for the operator $%
{\bf H}^{(0)}\boxplus \big( \bigoplus_{m=1}^nH_{\mathrm sym}^{(n,m)}\big)$ in the space $%
L_{\pi }^2(\Gamma _X)\otimes \big[
\bigoplus_{m=1}^nL^2\Psi_{\mathrm sym} ^n(X^m)\big]$. Thus, from
(\ref{mn}) we deduce the remaining statements of the theorem.\quad
$\blacksquare $\vspace{2mm}

\subsection{Harmonic forms}

In this section, we study the spaces ${\bf
K}^{(n)}:=\operatorname{Ker}{\bf H}^{(n)}$ of harmonic forms over
$\Gamma_X$. We give their description in terms of the spaces of
harmonic forms of the underlying manifold $X$. For this, we need
some auxiliary facts concerning Hilbert tensor algebras with
certain commutation relations.

\subsubsection*{Some   Hilbert tensor algebras}


Let ${\cal A}({\cal H}_1,\dots,{\cal H}_l)$ be the free Hilbert
tensor algebra generated by real separable  Hilbert spaces ${\cal
H}_1,\dots,{\cal H }_l$, $l\in\N$. That is, \begin{align*}{\cal
A}({\cal H}_1,\dots,{\cal H}_l)&:=\bigoplus_{m=0}^\infty {\cal
A}_m({\cal H}_1,\dots,{\cal H}_l) ,\\ {\cal A}_0({\cal
H}_1,\dots,{\cal H}_l)&:=\R,\\ {\cal A}_m({\cal H}_1,\dots,{\cal
H}_l)&:=\bigoplus_{i_1,\dots,i_m\in\{1,\dots,l\}}{\cal
H}_{i_1}\otimes\dots\otimes {\cal H}_{i_m},\qquad
m\in\N,\end{align*} with the usual addition and tensor product of
elements.

To each space ${\cal H}_i$, $i=1,\dots,l$, we associate a
parameter $p(i)\equiv p({\cal H}_i)\in\N$ (degree). Let $\Theta$
be the closure of the ideal in ${\cal A}({\cal H}_1,\dots,{\cal
H}_l)$ generated by the elements $$h\otimes
f-(-1)^{p(i)p(j)}f\otimes h,\qquad h\in{\cal H}_i,\, f\in{\cal
H}_j,\ i,j\in\{1,\dots,l\}.$$ That is,
\begin{multline*}\Theta:=\operatorname{c.l.s.}\big\{\,
a\otimes[h\otimes f-(-1)^{p(i)p(j)}f\otimes h]\otimes b\mid
\\ a, b\in{\cal A} ({\cal H}_1,\dots,{\cal H}_l),\
h\in{\cal H}_i,\, f\in{\cal H}_j,\
i,j\in\{1,\dots,l\}\,\big\},\end{multline*} where
$\operatorname{c.l.s.}$ means the closed linear span.

Let us define the quotient Hilbert  space $${\cal A}_{\mathrm
sym}({\cal H}_1,\dots,{\cal H}_l):={\cal A}({\cal H}_1,\dots,{\cal
H}_l)/\Theta.$$ As usual, we can identify ${\cal A}_{\mathrm
sym}({\cal H}_1,\dots,{\cal H}_l)$ with the orthogonal complement
of $\Theta$ in \linebreak ${\cal A}({\cal H}_1,\dots,{\cal H}_l)$.

\newcommand{\sign}{\operatorname{sign}}

\begin{lemma}\label{lemmas1} Let the linear continuous operator $\bf
P$ in ${\cal A}({\cal H}_1,\dots,{\cal H}_l)$ be defined through
the relation\begin{equation}\label{ersrwes} \begin{gathered} {\bf
P}(h_1\otimes\dots\otimes h_m):=\frac1{m!}\,\sum_{\sigma\in
S_m}\sign (\sigma,i_1,\dots,i_m)h_{\sigma(1)}\otimes\dots\otimes
h_{\sigma(m)},\\  h_k\in{\cal H}_{i_k},\ k=1,\dots,m,\
i_1,\dots,i_m\in\{1,\dots,l\}.\end{gathered}\end{equation}
Here\rom, $$\sign (\sigma,i_1,\dots,i_m):=\prod_{k<r:\,
\sigma(k)>\sigma(r)}(-1)^{p(i_{\sigma(k)})p(i_{\sigma(r)})},$$
with $\prod_{x\in\varnothing}a_x:=1$\rom. Then\rom, $\bf P$ is the
orthogonal projection of ${\cal A}({\cal H} _1,\dots,{\cal H}_l)$
onto \linebreak ${\cal A}_{\mathrm sym}({\cal H} _1,\dots,{\cal
H}_l)$\rom.
\end{lemma}

\noindent {\it Proof}. See
Appendix~B.\quad$\blacksquare$\vspace{2mm}

Let $$\Theta_m:=\Theta\cap {\cal A}_m({\cal H} _1,\dots,{\cal
H}_l)$$ and $$ {\cal A}_{m,\,{\mathrm sym}}({\cal H}
_1,\dots,{\cal H}_l):={\cal A}_m({\cal H }_1,\dots,{\cal
H}_l)/\Theta_m.$$ Evidently, $${\cal A}_{\mathrm sym }({\cal H}
_1,\dots,{\cal H}_l)=\bigoplus_{m=0}^\infty  {\cal
A}_{m,\,{\mathrm sym}}({\cal H }_1,\dots,{\cal H}_l).$$ The
following lemma gives an isomorphic description of the spaces
${\cal A}_{m,\,{\mathrm sym }}({\cal H} _1,\dots,{\cal H}_l)$.

\begin{lemma}\label{lemmas2} For each $m\in\N$\rom, there exists a
unitary isomorphism $${\cal U}_m:{\cal A}_{m,\,{\mathrm sym
}}({\cal H} _1,\dots,{\cal H}_l)\to \underset{\begin{gathered}
\scriptstyle s_1,\dots,s_l\in\Z_+\\ \scriptstyle
s_1+\dots+s_l=m\end{gathered}}{\bigoplus} \bigotimes _{i=1}^l
{\cal H}_i ^{\overset{p(i)}{\diamond}s_i}.$$ Here\rom, for each
$i\in\{1,\dots,l\}$\rom,  $\overset{p(i)}{\diamond}$ denotes the
antisymmetric tensor product $\wedge$ if $p(i)$ is odd and the
symmetric tensor product $\widehat\otimes$ if $p(i)$ is even\rom.
The unitary  operator ${\cal U}_m$ is constructed through the
relation
\begin{gather} {\cal U}_m\big({\bf
P}(f^{(1)}_1\otimes\dots\otimes f^{(1)}_{r_1}\otimes\dots\otimes
f^{(l)}_1\otimes\dots\otimes f^{(l)}_{r_l})\big):=\notag\\
=\bigg(\frac{m!}{r_1!\dotsm
r_l!}\bigg)^{1/2}\,\big(f^{(1)}_1\overset{p(1)}{\diamond}\dotsm
\overset{p(1)}{\diamond}f^{(1)}_{r_1}\big)
\otimes\dots\otimes\big(f^{(l)}_{1}\overset{p(l)}{\diamond}\dotsm
\overset{p(l)}{\diamond}f^{(l)}_{r_l}\big),\label{rtij}\\
f^{(i)}_k\in{\cal H}_i,\ k=1,\dots,r_i,\ r_1,\dots,r_l\in\Z_+,\
r_1+\dots+r_l=m, \notag\end{gather} the resulting operator ${\cal
U}_m$ being independent of the representation of a vector
from\linebreak  ${\cal A}({\cal H}_1,\dots,{\cal H}_l)$\rom.
\end{lemma}

\noindent {\it Proof}.  See Appendix B. \quad $\blacksquare$
\vspace{2mm}

Now, for each $n\in\N$, we define the subspace ${\cal A}^{n}({\cal
H}_1,\dots,{\cal H}_l)$ of ${\cal A}({\cal H}_1,\dots,{\cal H}_l)$
by setting \begin{align*}{\cal A}^{n}({\cal H}_1,\dots,{\cal
H}_l)&:=\bigoplus_{m=1}^n {\cal A}^{n}_m({\cal H}_1,\dots,{\cal
H}_l),\\ {\cal A}^{n}_m({\cal H}_1,\dots,{\cal H}_l)&:=
\underset{\begin{gathered} \scriptstyle
i_1,\dots,i_m\in\{1,\dots,l\}\\ \scriptstyle
p(i_1)+\dots+p(i_m)=n\end{gathered}}{\bigoplus} {\cal
H}_{i_1}\otimes \dots\otimes {\cal H}_{i_m}.\end{align*} Let also
\begin{align*} \Theta^n&:=\Theta\cap{\cal A}^n({\cal H}_1,\dots,{\cal
H}_l),\\ {\cal A}^n_{\mathrm sym}({\cal H}_1,\dots,{\cal
H}_l)&:={\cal A}^n ({\cal H}_1,\dots,{\cal
H}_l)/\Theta^n.\end{align*} Evidently, \begin{align}{\cal
A}^n_{\mathrm sym}({\cal H}_1,\dots,{\cal H}_l)&=\bigoplus_{m=1}^n
{\cal A}^n_{m,\,{\mathrm sym}}({\cal H}_1,\dots,{\cal
H}_l),\label{rdtuztzdd}\\ {\cal A}^n_{m,\,{\mathrm sym}}({\cal
H}_1,\dots,{\cal H}_l):&={\cal A}^n_{m}({\cal H}_1,\dots,{\cal
H}_l)/\Theta_m^n,\qquad \Theta_m^n:=\Theta\cap  {\cal
A}^n_{m}({\cal H}_1,\dots,{\cal H}_l).\notag\end{align} By
Lemma~\ref{lemmas1}, the orthogonal  projection $ {\bf P}^n_m$ of
${\cal A}^n_{m}({\cal H}_1,\dots,{\cal H}_l)$ onto ${\cal
A}^n_{m,\,{\mathrm sym}}({\cal H}_1,\dots,{\cal H}_l)$ is the
restriction of ${\bf P}$ to ${\cal A}^n_{m}({\cal H}_1,\dots,{\cal
H}_l)$, and by \eqref{rdtuztzdd} and Lemma~\ref{lemmas2}  the
restrictions of the  ${\cal U}_m$'s, $m=1,\dots,n$, define the
unitary operator
\begin{equation}\label{qawareaa} {\cal U}^n:{\cal A}^n_{{\mathrm sym
}}({\cal H} _1,\dots,{\cal H}_l)\to
\bigoplus_{m=1}^n\underset{\begin{gathered} \scriptstyle
s_1,\dots,s_l\in\Z_+\\ \scriptstyle s_1+\dots+s_l=m\\ \scriptstyle
p(1)s_1+\dots+p(l)s_l=n\end{gathered}}{\bigoplus} \bigotimes
_{i=1}^l {\cal H}_i ^{\overset{p(i)}{\diamond}s_i}.
\end{equation}

\begin{remark}\rom{ Actually, ${\cal U}^n$ is a natural
isomorphism generated by the passage to summation in ordered
families of indices in the definition of ${\cal A}^n_{\mathrm
sym}({\cal H}_1,\dots,{\cal H}_l)$, which uses the commutation
relation $$h\otimes f=(-1)^{p(i)p(j)}f\otimes h,\qquad h\in{\cal
H}_i,\ f\in{\cal H}_j,\ i,j\in\{1,\dots,l\}.$$ }\end{remark}

\begin{remark}\label{drtrzdtftz}\rom{Setting ${\cal A}^0_{\mathrm
sym}({\cal H}_1,\dots,{\cal H}_l):=\R$, one gets the orthogonal
decomposition $${\cal A}_{\mathrm sym}({\cal H}_1,\dots,{\cal
H}_l)=\bigoplus_{n=0}^\infty{\cal A}^{n}_ {\mathrm sym}({\cal
H}_1,\dots,{\cal H}_l).$$}\end{remark}

\subsubsection*{The Kernel of the Hodge--deRham Laplacian}

Our next goal is to investigate the kernel of $ {\bf H}^{(n)}$. We
first need  the following general result.

\begin{lemma} \label{awsedrgvgvdfgrdrtf} Let $A$ and $B$ be
self-adjoint\rom, non-negative operators in separable Hilbert
spaces $\cal H$ and $\cal K$\rom, respectively\rom. Then\rom, we
have $$\Ker (A\boxplus B)=\Ker A\otimes \Ker B,$$ where $A\boxplus
B$ is the closure of the operator $A\otimes I+ I\otimes B$ from
the set $\operatorname{Dom}A\otimes_{\mathrm
a}\operatorname{Dom}B$\rom.\end{lemma}

\noindent{\it Proof}.  $\Ker A$ and $\Ker B$ are closed subspaces
of $\cal H$, resp.\ $\cal K$, and so  their tensor product $\Ker
A\otimes \Ker B$ is a closed subspace of the space ${\cal
H}\otimes{\cal K}$. The inclusion $\Ker A\otimes\Ker B\subset \Ker
(A\boxplus B)$ is trivial. Let $f\in \Ker (A\boxplus B)$. Using
the theory of operators admitting separation of variables (e.g.\
\cite[Ch.~6]{B}), we have
\begin{align}0=(A\boxplus B
f,f)&=\int_{\R_+^2}(x_1+x_2)\,d(E(x_1,x_2)f,f)\notag\\ &=\int_{\R^
2_+}x_1\,d(E(x_1,x_2)f,f)+\int_{\R^
2_+}x_2\,d(E(x_1,x_2)f,f)\notag\\ &=(A\otimes If,f)+(I\otimes
Bf,f),\label{dderuz}\end{align} where $E$ is the joint resolution
of the identity of the commuting operators $A\otimes I$ and
$I\otimes B$. Since both operators $A\otimes I$ and $I\otimes B$
are non-negative, we conclude from \eqref{dderuz} that $$f\in
\Ker( A\otimes I)\cap\Ker (I\otimes B)=\Ker A\otimes\Ker
B.\quad\blacksquare\vspace{2mm}$$

Let us fix any $i_1,\dots,i_m\in\{1,\dots,d\}$, $i_1+\dots+i_m=n$.
For any $\om_r\in\Omega^{i_r}(X)$, $r=1,\dots,m$, we define the
form  $$X^m \ni(x_1,\dots,x_m)\mapsto \tilde
\omega_r(x_1,\dots,x_m):=\omega_r(x_r)\in\wedge^{i_r}(T_{x_r}X)\subset\wedge
^{i_r}(T_{(x_1,\dots,x_m)}X^m).$$ Now, we set $$U_{i_1,\dots,i_m}
(\omega_1\otimes\dots\otimes\omega_m):=\bigg(\frac{n!}{i_1!\dotsm
i_m! }\bigg)^{1/2}\,\tilde \omega_1\wedge\dots\wedge
\tilde\omega_m\in\Omega^n(X^m).$$ (We use here the convention that
the exterior product of two forms, $\omega$ and $\nu$, is given by
$\omega\wedge\nu:=\operatorname{AS}(\omega\otimes\nu)$, where
$\operatorname{AS}$ denotes the antisymmetrization operator). It
is easy to see that $U_{i_1,\dots,i_m}$ can be extended by
linearity and continuity to a linear isometric operator
$$U_{i_1,\dots,i_m}: L^2\Omega^{i_1}(X)\otimes\dots\otimes
L^2\Omega^{i_m}(X)\to L^2\Omega^n(X^m)$$ with the image
$$\operatorname{Im}U_{i_1,\dots,i_m}=
L^2\Psi_{i_1,\dots,i_m}(X^m),$$ where
$L^2\Psi_{i_1,\dots,i_m}(X^m)$ denotes the space of the  forms $$
X^m\ni(x_1,\dots,x_m)\mapsto \om(x_1,\dots,x_m)\in
(T_{x_1}X)^{\wedge i_1}\wedge\dots\wedge (T_{x_m}X)^{\wedge i_m}$$
that are square integrable with respect to $dx_1\dotsm dx_m$.

Setting \begin{equation} L^2\Psi^n(X^m):=
\underset{\begin{gathered} \scriptstyle
i_1,\dots,i_m\in\{1,\dots,d\}\\ \scriptstyle
i_1+\dots+i_m=n\end{gathered}}{\bigoplus}
L^2\Psi_{i_1,\dots,i_m}(X^m),\label{drzgu}\end{equation} we
construct,  by using the $U_{i_1,\dots,i_m}$'s, the unitary
isomorphism $$U_m^n:\underset{\begin{gathered} \scriptstyle
i_1,\dots,i_m\in\{1,\dots,d\}\\ \scriptstyle
i_1+\dots+i_m=n\end{gathered}}{\bigoplus}
L^2\Omega^{i_1}(X)\otimes \dots\otimes L^2\Omega^{i_m}(X)\to
L^2\Psi^n(X^m),$$
or equivalently
\begin{equation}\label{esres}U_m^n : {\cal
A}^n_m(L^2\Omega^1(X),\dots,L^2\Omega^d(X))\to
L^2\Psi^n(X^m),\end{equation} where $p(i)=p(L^2\Omega^i(X)):=i$.

We notice that the restriction of the orthogonal projection
$$P_m^n: L^2\Omega^n(X^m)\to L^2\Omega^n_{\mathrm sym}(X^m)$$ to
the subspace $L^2\Psi^n(X^m)$ determines the orthogonal projection
$$P_m^n: L^2\Psi^n(X^m)\to L^2\Psi^n_{\mathrm sym}(X^m).$$

\begin{lemma} \label{lemmas3} We have $$P_m^nU_m^n=U_m^n{\bf P}_m^n,$$
where ${\bf P}_m^n$ is the  orthogonal projection of ${\cal
A}_m^n(L^2\Omega^1(X),\dots,L^2\Omega^d(X))$ onto \linebreak
${\cal A}^n_{m,\mathrm
sym}(L^2\Omega^1(X),\dots,L^2\Omega^d(X))$\rom.
\end{lemma}

\noindent {\it Proof}. For any $\om_r\in\Omega^{i_r}(X)$,
$r=1,\dots,m$, $i_1,\dots,i_r\in\{1,\dots,d\}$, $i_1+\dots+i_m=n$,
we get by using Lemma~\ref{lemmas1}
\begin{align*}&(P_m^nU_m^n\om_1\otimes\dots\otimes
\om_m)(x_1,\dots,x_m)=\\ &\qquad=\bigg(\frac {n!}{i_1!\dotsm
i_m!}\bigg)^{1/2}\sum_{\sigma\in
S_m}\om_1(x_{\sigma(1)})\wedge\dotsm\wedge
\om_m(x_{\sigma(m)})\\&\qquad =\bigg(\frac {n!}{i_1!\dotsm
i_m!}\bigg)^{1/2}\sum_{\sigma\in S_m}
\sign(\sigma,i_1,\dots,i_m)\,\om_{\sigma(1)}(x_1)\wedge\dotsm\wedge
\om_{\sigma(m)}(x_m)\\&\qquad =\sum_{\sigma\in
S_m}\sign(\sigma,i_1,\dots,i_m)\,(U_m^n
\om_{\sigma(1)}\otimes\dotsm\otimes
\om_{\sigma(m)})(x_1,\dots,x_m)\\&\qquad =(U_m^n{\bf
P}_m^n\om_1\otimes \dots\otimes \om_m)(x_1,\dots,x_m).\qquad
\blacksquare\end{align*}

Since $P^n_m$ is the orthogonal projection of $L^2\Psi^n(X^m)$
onto $L^2\Psi^n_{\mathrm sym}(X^m)$ and ${\bf P}_m^n$ is the
orthogonal projection of ${\cal
A}_m^n(L^2\Omega^1(X),\dots,L^2\Omega^d(X))$ onto ${\cal
A}_{m,\,\mathrm sym}^n(L^2\Omega^1(X),\dots,L^2\Omega^d(X))$, we
conclude from \eqref{esres} and Lemma~\ref{lemmas3} that the
restriction of $U_m^n$ to  \linebreak ${\cal A}_{m,\,\mathrm
sym}^n(L^2\Omega^1(X),\dots,L^2\Omega^d(X))$ defines the unitary
isomorphism $$U_m^n:  {\cal A}_{m,\,\mathrm
sym}^n(L^2\Omega^1(X),\dots,L^2\Omega^d(X))\to L^2\Psi^n_{\mathrm
sym }(X^m).$$

Finally, setting
\begin{equation}\label{uztftf}U^n:=\bigoplus_{m=1}^n
U^n_m,\end{equation} we get the unitary mapping $$U^n:  {\cal
A}_{\mathrm sym}^n(L^2\Omega^1(X),\dots,L^2\Omega^d(X))\to
\bigoplus_{m=1}^n L^2\Psi^n_{\mathrm sym }(X^m).$$

We denote by ${\cal K}^{(i)}$ the kernel of the Hodge--deRham
Laplacian $H^{(i)}$ in the space $L^2\Omega^i(X)$, $i=1,\dots,d$.
Each ${\cal K}^{(i)}$ as a closed subspace of the Hilbert space
$L^2\Omega^i(X)$ is itself a Hilbert space. Let also  ${\bf
K}^{(n)}$ denote the kernel of the operator ${\bf H}^{(n)}$.

\begin{theorem}\label{tfgzitzz} We have \begin{equation}\label{uiiiuoio}
 I^{n}{\bf
K}^{(n)} =\{\operatorname{const}\}\otimes \left[ U^n {\cal A
}^n_{\mathrm sym}({\cal K}^{(1)},\dots,{\cal
K}^{(d)})\right],\end{equation} where $p(i):=i$\rom,
$i=1,\dots,d$\rom.
\end{theorem}

\noindent {\it Proof}. By \cite[Theorem~4.3]{AKR1},
\begin{equation}\label{rtdttz}{\bf K}^{(0)}:=\operatorname{Ker}{\bf H}^{(0)}=\{\operatorname{const}\},
\end{equation}
and hence by Theorem~\ref{thonsa2} and
Lemma~\ref{awsedrgvgvdfgrdrtf}
\begin{equation}\label{einss}I^{n}{\bf
K}^{(n)}=\{\operatorname{const}\}\otimes\left[\bigoplus_{m=1}^n\operatorname{Ker}
H^{(n,m)}_{\mathrm sym}\right].\end{equation}

Let us find the kernel of the Hodge--deRham Laplacian  $H^{(n,m)}$
in the space $L^2\Psi^n(X^m)$. The operator $H^{(n,m)}$ acts
invariantly in each space in the direct sum \eqref{drzgu}, so that
it suffices to find the kernel of each restriction
$H^{(n,m)}_{i_1,\dots,i_m}$ of $H^{(n,m)}$ to the subspace
$L^2\Psi_{i_1,\dots,i_m}(X^m)$.

By using the operator $U_m^n$, we easily conclude that
$$(U_m^n)^{-1}H^{(n,m)}_{i_1,\dots,i_m}U_m^n=(c_1H^{(i_1)})\boxplus\dots\boxplus
(c_m H^{(i_m)}),$$ where $c_1,\dots,c_m$ are non-zero constants.
Therefore, by Lemma~\ref{awsedrgvgvdfgrdrtf} $$\operatorname{Ker}
H^{(n,m)}_{i_1,\dots,i_m}=U_m^n({\cal K}^{(i_1)}\otimes
\dots\otimes {\cal K}^{(i_m)}),$$ which yields that
\begin{align*}\operatorname{Ker} H^{(n,m)}&=U_m^n
\underset{\begin{gathered} \scriptstyle
i_1,\dots,i_m\in\{1,\dots,d\}\\ \scriptstyle
i_1+\dots+i_m=n\end{gathered}}{\bigoplus} ({\cal K}^{(i_1)}\otimes
\dots\otimes {\cal K}^{(i_m)})\\&=U_m^n{\cal A} _m^n({\cal
K}^{(1)},\dots,{\cal K}^{(d)}).\end{align*} Since
$$H^{(n,m)}_{\mathrm sym}P_m^n=P_m^n H^{(n,m)},$$ we get
$$\operatorname{Ker}H^{(n,m)}_{\mathrm sym}=P_m^{(n)}
\operatorname{Ker}H^{(n,m)},$$ which implies by
Lemma~\ref{lemmas3} that
\begin{equation}\operatorname{Ker}H^{(n,m)}_{\mathrm sym}=U_m^n
{\cal A}^n_{m,\,\mathrm sym}({\cal K}^{(1)},\dots,{\cal
K}^{(d)}).\label{zweii}\end{equation} Combining \eqref{uztftf},
\eqref{einss} and \eqref{zweii}, we get the conclusion of the
theorem.\qquad $\blacksquare$

\begin{corollary} The isomorphisms $I^n$\rom, $U^n$ and the
equality \eqref{rtdttz} generate the unitary isomorphism  of the
Hilbert spaces $$\bigoplus_{n=0}^\infty {\bf K}^{(n)}\simeq{\cal
A}_{\mathrm sym}({\cal K}^{(1)},\dots,{\cal
K}^{(d)}).$$\end{corollary}

\noindent {\it Proof}. For each $n\in\N$, we get from
\eqref{uiiiuoio} the unitary isomorphism of the spaces $${\bf
K}^{(n)}\simeq{\cal A}^n_{\mathrm sym}({\cal K}^{(1)},\dots,{\cal
K}^{(d)}).$$ Moreover, it follows from \eqref{rtdttz}  that ${\bf
K }^{(0)}\simeq\R$. Hence, the conclusion of the corollary follows
from Remark~\ref{drtrzdtftz}.\quad$\blacksquare$

\begin{remark}\rom{ Formula
\eqref{uiiiuoio} is wrong in the case where the manifold $X$ has
finite volume (in that case the Poisson measure $\pi$ is
concentrated on the space of finite configurations over $X$).
Instead of \eqref{uiiiuoio}, one then gets $$ I^{n}{\bf K}^{(n)}
=\operatorname{Ker}{\bf H }^{(0)}\otimes \left[ U^n {\cal A
}^n_{\mathrm sym}({\cal K}^{(1)},\dots,{\cal K}^{(d)})\right],$$
the space $\operatorname{Ker}{\bf H} ^{(0)}$ being
infinite-dimensional. }\end{remark}

\subsection{Structure of $L^2$-cohomologies}

The aim of this section is to study the structure of the spaces ${\cal H}%
_\pi ^n$ of $L^2$-cohomologies of $\Gamma _X$ using the
representation of the kernel of ${\bf H}^{(n)}$ given by
Theorem~\ref{tfgzitzz}. The following proposition reflects a quite
standard fact in the $L^2$-theory.

\begin{proposition}
The natural isomorphism between ${\cal H}^n_\pi$ and the
orthogonal complement of $B^n$ to $Z^n$ is the isomorphism of the
Hilbert spaces
\begin{equation}
{\cal H}_\pi ^n\simeq \operatorname{Ker}{\bf H}^{(n)}.
\label{harm1}
\end{equation}
\end{proposition}

\noindent {\it Proof}. Using \cite[Proposition~A.1]{Ar0}, we
conclude from Proposition~\ref{guzftdrt} and formula
\eqref{pserdo} that
\begin{equation}
L_\pi ^2\Omega ^n= \operatorname{Ker}{\bf H}^{(n)}\oplus \overline{\operatorname{Im}
{\bf d}_{n-1}}%
\oplus \overline{\operatorname{Im}{\bf d}_n^{*}}  \label{hdr}
\end{equation}
(weak Hodge--deRham decomposition). For the closed operator $\bar
{\bf d}_n$ we have the standard decomposition
\[
L_\pi ^2\Omega ^n=\operatorname{Ker}\bar {\bf d}_n\oplus
\overline{\operatorname{Im} {\bf d}_n^{*}},
\]
which together with (\ref{hdr}) implies the result. $\blacksquare
$\vspace{2mm}

Due to the Hodge--deRham theory of the underlying manifold $X$, we
have the isomorphisms
\begin{equation}\label{rrttz}
{\cal K}^k\simeq {\cal H}^k_{(2)}(X),\qquad k=1,\dots,d,
\end{equation}
where ${\cal H}^k_{(2)}(X){:=}\operatorname{Ker}
d_k/\,\overline{\operatorname{Im}d_{k-1}}$ ($d_j$, $j=1,\dots,d$,
denoting the Hodge differential of $X$) is the corresponding space
of (reduced) $L^2$-cohomologies of $X$.

\begin{remark}\rom{Because of the elliptic regularity of the
Hodge--deRham Laplacian on $X$, there exists a canonical map
${\cal H}^*_{(2)}(X)\to {\cal H}^*(X)$, where ${\cal H}^*(X)$ is
the deRham cohomology of $X$.  In general, this map is neither
surjective, nor injective. }\end{remark}

\begin{theorem}
\label{betti} \rom{1)} The isomorphisms \eqref{harm1}\rom,
$I^n$\rom, $U^n$\rom, ${\cal U}^n$\rom, and \eqref{rrttz} generate
the unitary isomorphism of the Hilbert spaces
\begin{equation}\label{cohom} {\cal H}^n_\pi\simeq \bigoplus_{m=1}^n
\bigoplus _{1\le k_1<\dots< k_m\le d }\underset{\begin{gathered}
\scriptstyle s_1,\dots,s_m\in\N\\
 \scriptstyle
k_1 s_1+\dots+k_m s_m=n\end{gathered}}{\bigoplus} ({\cal
H}_{(2)}^{k_1}(X))^{\overset{k_1}{\diamond}s_1}\otimes\dots\otimes
({\cal
H}_{(2)}^{k_m}(X))^{\overset{k_m}{\diamond}s_m}.\end{equation}

\rom{2)} Let  $\beta_k{:=}\operatorname{dim} {\cal
H}^k_{(2)}(X)<\infty$\rom, $k=1,\dots,d$\rom. Then\rom, all the
spaces ${\cal H}_\pi^n$\rom, $n\in\N$\rom, are
finite-dimensional\rom, and we have the following formula for
their dimensions $b_n$\rom:
\begin{equation}\label{weerdt} b_n= \sum_{m=1}^n  \, \sum _{1\le k_1<\dots< k_m\le d
}\underset{\begin{gathered} \scriptstyle s_1,\dots,s_m\in\N\\
 \scriptstyle
k_1 s_1+\dots+k_m s_m=n\end{gathered}}{\sum} \beta
_{k_1}^{(s_1)}\dotsm\beta _{k_m}^{(s_m)},\end{equation} where
\begin{equation}\label{tfttz} \beta _k^{(s)}:=\begin{cases}
\binom{ \beta _k}{ s} ,& k=1,3,\dots,\\ \binom {\beta _k+s-1 }{ s
},& k=2,4,\dots\end{cases}\end{equation}
\end{theorem}

\noindent{\it Proof}. 1) Follows from Theorem~\ref{tfgzitzz}.
Actually, \eqref{cohom} is a more explicit form of \eqref{uztftf}.

2) It is easy to see that, for a finite-dimensional space ${\cal
H}$, we have
\begin{align*}
\operatorname{dim}{\cal H}^{\widehat{\otimes }s}
&=\frac{(s+1)(s+2)\dotsm
(s+\operatorname{dim}{\cal H}-1)}{%
(\operatorname{dim}{\cal H}-1)!}=\binom{ \operatorname{dim}{\cal
H}+s-1}  s , \\   \operatorname{dim}{\cal H}^{\wedge s} &=\binom
{\operatorname{dim}{\cal H}} s .\end{align*} The statement follows
now from (\ref{cohom}). \quad $\blacksquare
$

\begin{corollary}\label{qwqqwqwqq}
Let  $\beta _1,\dots,\beta _d$  be finite\rom, and moreover let
$\beta _k=0$ for all $k$ even\rom. Then\rom:
\begin{align*}
b_k &=0,\qquad \text{\rom{for all} } k>K_0:=\sum_{i=1}^{d}i\beta
_i,
\\ b_{K_0} &=1.  \end{align*}
\end{corollary}

\noindent {\it Proof}. The condition $\beta_k=0$ for all $k$ even
implies that $$ {\cal H}^n_\pi\simeq \bigoplus_{m=1}^n
\underset{\begin{gathered}\scriptstyle 1\le k_1<\dots< k_m\le d \\
\scriptstyle  k_1,\dots,k_m\ {\mathrm odd}
\end{gathered}}{\bigoplus}\,
\underset{\begin{gathered}
\scriptstyle s_1,\dots,s_m\in\N\\
 \scriptstyle
k_1 s_1+\dots+k_m s_m=n\end{gathered}}{\bigoplus} ({\cal
H}_{(2)}^{k_1}(X))^{\wedge s_1}\otimes\dots\otimes ({\cal
H}_{(2)}^{k_m}(X))^{\wedge s_m}.$$ Obviously $( {\cal
H}_{(2)}^k(X))
 ^{\wedge s}=0$ for $s>\beta _k$ and $%
( {\cal H}^k(X)) ^{\wedge s}= {\Bbb R}^1$ for $s=\beta _k$, which
implies the result.\quad $\blacksquare$

\begin{example}\label{sesesemn}\rom{ Let $X$ be a manifold with a cylindrical end
(that is, $X=M\cup (N\times \R_+^1)$ for some compact manifold $M$
with boundary $N$). It is proven in \cite{APS} that  ${\cal
H}^k_{(2)}(X)$ is isomorphic to the image of the canonical map
${\cal H}_0^k(X)\to {\cal H}^k(X)$, where ${\cal H}_0^k(X)$ is the
space of the compactly supported deRham cohomologies of $X$.  By
e.g.\ \cite{Bott}, the spaces ${\cal H}^k(X)$ are
finite-dimensional. Thus, all ${\cal H}_{(2)}^k(X)$ are
finite-dimensional and, in general, non-trivial, and hence so are
all  spaces  ${\cal H}_\pi^n$. For a bigger class of examples of
manifolds $X$ with finite-dimensional spaces ${\cal H}^k_{(2)}(X)$
see \cite{Mul}. }\end{example}

\begin{example}\rom{ Let $d=2$. Then, $\beta_0=\beta_2=0$ (see e.g.\
\cite{And}), and if $X$ is as in Example~\ref{sesesemn}, we also
have $\beta_1<\infty$. Thus,  $X$ satisfies the conditions of
Corollary~\ref{qwqqwqwqq}, and we have $b_k=0$ for all $k>\beta_1$
and $b_k=\binom{\beta_1}{k}$ for $k\le \beta_1$.}\end{example}

\begin{remark}\rom{ The vanishing of the spaces ${\cal H}_\pi^n$ does
not, in general, imply the absence of non-exact closed forms.
Suppose, for example, that $X=\R^1$. Clearly, there are no
$L^2$-harmonic forms on $\R^1$, which implies that all the spaces
${\cal H}^n_\pi (\Gamma_{\R^1})$ are trivial. Let us consider a
1-form $\varphi(x)=g(x)\, dx$ on $\R^1$ such that $g(x)$ has a
compact support  and $\int_{\R^1} \varphi\ne 0$. The latter
implies that $\varphi\ne d_0f$ for any $f\in L^2(\R^1)$. We now
define $\Phi\in L^2_\pi \Omega^1(\Gamma_{\R^1})$ setting
$\Phi(\gamma)_x{:=}\varphi(x)$. It is easy to see that $\Phi\ne
{\bf d}_0F$ for any $F\in L^2_\pi(\Gamma_{\R^1})$ and ${\bf
d}_1\Phi=0$.}\end{remark}

\begin{example}\rom{
{\bf Marked configuration spaces. }Let $Y=X\times M$, where $M$ is
a compact Riemannian manifold. We note that $\Gamma _Y$ coincides
up
to a set of zero $\pi$ measure with the marked configuration space $\Gamma _{%
X}(M)$, see e.g.\  \cite{kingman}. Let us recall that the latter
space is defined as follows: $$\Gamma_{X}(M):=\big\{\,
\gamma\in\Gamma_{X\times M}:\forall (x_1,m_1),(x_2,m_2)\in\gamma:\
(x_1,m_1)\ne(x_2,m_2)\Rightarrow x_1\ne x_2 \,\big\}. $$ The
K\"{u}nneth formula  implies
\begin{align*}
{\cal H}_{(2)}^n( X\times M) &=\bigoplus_{m=0}^n {\cal
H}_{(2)}^m(X)\otimes {\cal H}_{(2)}%
^{n-m}(M) . \end{align*} We remark that, for each $k$, ${\cal
H}^k_{(2)}(M)={\cal H}^k(M)$ and is finite-dimensional. Thus, all
the spaces ${\cal H}_\pi^{n}(\Gamma_X(M))$ are finite-dimensional,
provided so are all ${\cal H}^k_{(2)}(X)$. }
\end{example}

\section{Appendix}

\subsection{Appendix~A: Laplacian on the configuration space}

We recall here the definition of the Laplacian on the
configuration space and some facts about it from \cite{AKR1},
which we present in a form adapted to the aims of the present
paper (see also \cite{ADL1,ADL2}).

Let $F:\Gamma_X\to\R$. For fixed $\gamma\in\Gamma_X$ and
$x\in\gamma$, we define the function $${\cal O}_{\gamma,x}\ni
y\mapsto F_x(\gamma,y):=F((\gamma\setminus\{x\})\cup\{y\})\in\R.$$
We say that $F$ is differentiable at $\gamma\in\Gamma_X$ if, for
each $x\in\gamma$, the function $F_x(\gamma,\cdot)$ is
differentiable at $x$ and $$\nabla^\Gamma F(\gamma):=(\nabla^X F_x
(\gamma,x))_{x\in\gamma}\in T_\gamma\Gamma_X. $$ Analogously, the
higher order derivatives of $F$ are defined,
$(\nabla^\Gamma)^{(m)}F(\gamma)\in(T_\gamma\Gamma_X)^{\otimes m}$,
$m\in\N$.

A function $F:\Gamma_X\to\R$ is called local if there exists a
compact $\Lambda\subset X$ such that $F(\gamma)=F(\gamma_\Lambda)$
for each $\gamma\in\Gamma_X$.

We define ${\cal FC}_{\mathrm b}^\infty({\cal D},\Gamma_X)$ as the
set of all functions $F:\Gamma_X\to\R$ of the form
\begin{equation}\label{errdtftfz}
F(\gamma)=g_F(\langle\varphi_1,\gamma\rangle,\dots,\langle\varphi_N,\gamma\rangle),
\end{equation}
where $g_F\in C^\infty_{\mathrm b}(\R^N)$ and
$\varphi_1,\dots,\varphi_N\in{\cal D}{:=} C_0^\infty(X)$(:$=$the
set of all infinitely differentiable functions on $X$ with compact
support). Each function $F\in {\cal FC}_{\mathrm b}^\infty({\cal
D},\Gamma_X)$ is evidently bounded, local, and infinitely
differentiable with derivatives satisfying the estimate
$$\|(\nabla^\Gamma)^{(m)}F(\gamma)\|_{(T_\gamma\Gamma_X)^{\otimes
m }} \le \langle \varphi^{\otimes m},\gamma^{\otimes
m}\rangle\qquad\text{for all }\gamma\in\Gamma_X,$$ with some
$\varphi\in C_0(X)$ depending on $F$ and $m\in\N$.

On the space $L^2_\pi(\Gamma_X)$ we consider the pre-Dirichlet
form $${\cal E}_\pi^{(0)}(F_1,F_2):=\int_{\Gamma_X}\langle
\nabla^\Gamma F(\gamma),\nabla^\Gamma
F(\gamma)\rangle_{T_\gamma\Gamma_X}\,\pi(d\gamma)$$ with domain
$\operatorname{Dom}{\cal E}_\pi^{(0)}{:=}{\cal FC}_{\mathrm
b}^\infty({\cal D},\Gamma_X)$, which is dense in
$L^2_\pi(\Gamma_X)$.

The following theorem can be  proved by using formula \eqref{3.1}.

\begin{theorem}\label{guzserser}For any $F_1,F_2\in {\cal FC}_{\mathrm b}^\infty({\cal
D},\Gamma_X)$\rom, we have $${\cal
E}^{(0)}_\pi(F_1,F_2)=\int_{\Gamma_X}({\bf H}^{(0)}F_1)(\gamma)
F_2(\gamma)\,\pi(d\gamma).$$ Here\rom, ${\bf
H}^{(0)}=-\Delta^\Gamma$ is the operator in $L^2_\pi(\Gamma_X)$
with domain $\operatorname{Dom}{\bf H}^{(0)}{:=}{\cal FC}_{\mathrm
b}^\infty({\cal D},\Gamma_X)$ that is given by the formula
\begin{equation}\label{esrr}({\bf
H}^{(0)}F)(\gamma):=-\sum_{x\in\gamma}\Delta^X
F_x(\gamma,x),\qquad F\in{\cal FC},\end{equation} $\Delta^X$
denoting the Laplacian on $X$\rom.\end{theorem}

From Theorem~\ref{guzserser} we conclude that the bilinear form ${\cal E}%
_\pi ^{(0)}$ is closable in the space $L_\pi ^2(\Gamma_X)$. The
generator of its
closure (being actually the Friedrichs extension of the operator ${\bf H}%
^{(0)}$, for which we preserve the same notation) will be called
the Laplacian on $\Gamma _X$. By \eqref{esrr}, ${\bf H}^{(0)}$ is
the lifting of the  Laplacian on $X$.

\begin{theorem} The operator ${\bf H}^{(0)}$ is essentially
self-adjoint on ${\cal FC}_{\mathrm b}^\infty({\cal
D},\Gamma_X)$\rom.\end{theorem}

\noindent {\it Proof}. See \cite[Theorem~5.3]{AKR1}.\quad
$\blacksquare$

\subsection{Appendix B: Proof of Lemmas \ref{lemmas1} and
\ref{lemmas2}}

We first prove Lemma \ref{lemmas1}.  Extending the
relation~\eqref{ersrwes} by linearity and continuity, we get a
linear continuous operator $\bf P$ in ${\cal A}({\cal
H}_1,\dots,{\cal H}_l)$.

Let us show that the operator $\bf P$ is self-adjoint. For
arbitrary $f_k\in{\cal H}_{i_k}$ and $g_k\in{\cal H}_{j_k}$,
$i_k,j_k\in \{1,\dots,l\}$, $k=1,\dots,m$,  we get from
\eqref{ersrwes}:
\begin{gather*}({\bf P}(f_1\otimes\dots\otimes f_m),
g_1\otimes\dots\otimes g_m)=\frac 1{m!}\,\sum_{\sigma\in S_m}
\sign(\sigma,i_1,\dots,i_m)\prod_{k=1}^m (f_{\sigma(k)},g_k)\\
=\frac 1{m!}\, \sum_{\sigma\in S_m}
\sign(\sigma^{-1},i_1,\dots,i_m)\prod_{k=1}^m
(f_k,g_{\sigma(k)})\\
 = \frac{1}{m!}\,
\sum_{\sigma\in S_m:\; i_1=j_{\sigma(1)},\dots,
i_m=j_{\sigma(m)}}\;\prod_{r<s:\
\sigma^{-1}(r)>\sigma^{-1}(s)}(-1)^{p(i_{\sigma^{-1}(r)})p(i_{\sigma^{-1}(s)})}
\,\prod_{k=1}^m (f_k,g_{\sigma(k)})\\ =\frac{1}{m!}\,
\sum_{\sigma\in S_m:\; i_1=j_{\sigma(1)},\dots,
i_m=j_{\sigma(m)}}\;\prod_{r<s:\ \sigma^{-1}(r)>\sigma^{-1}(s)}
(-1)^{p(j_r)p(j_s)}\,\prod_{k=1}^m (f_k,g_{\sigma(k)})\\
=\frac{1}{m!}\, \sum_{\sigma\in S_m:\; i_1=j_{\sigma(1)},\dots,
i_m=j_{\sigma(m)}}\;\prod_{r<s:\ \sigma(r)>\sigma(s)}
(-1)^{p(j_{\sigma(r)})p(j_{\sigma(s)})}\,\prod_{k=1}^m
(f_k,g_{\sigma(k)})
\\ =(f_1\otimes
\dots\otimes f_m,{\bf P}(g_1\otimes\dots\otimes g_m) ),
\end{gather*} and so $\bf P$ is indeed self-adjoint.

Next, it follows from the definition of $\bf P$ that, for
$f_k\in{\cal H}_{i_k}$, $i_k\in\{1,\dots,l\}$, $k=1,\dots,m$,
\begin{multline}\label{uuzfrdt}
{\bf P}(f_1\otimes \dots \otimes f_r\otimes
f_{r+1}\otimes\dots\otimes f_m)\\ = (-1)^{p(i_r)p(i_{r+1})}{\bf
P}(f_1\otimes\dots\otimes f_{r-1}\otimes f_{r+1}\otimes f_r\otimes
f_{r+2}\otimes \dots\otimes f_m).\end{multline} The latter formula
implies that, for each $\sigma\in S_m$,
\begin{equation}\label{ethess}{\bf P}(f_{\sigma(1)}\otimes\dots\otimes
f_{\sigma(m)})=\sign(\sigma,i_1,\dots,i_m){\bf
P}(f_1\otimes\dots\otimes f_m),\end{equation} and hence
\begin{align*} {\bf P}^2 (f_1\otimes\dots\otimes
f_m)&=\frac1{m!}\sum_{\sigma\in
S_m}\sign(\sigma,i_1,\dots,i_m){\bf
P}(f_{\sigma(1)}\otimes\dots\otimes f_{\sigma(m)})\\
&=\frac{1}{m!}\, \sum_{\sigma\in S_m
}\sign(\sigma,i_1,\dots,i_m)^2\,
 {\bf P} (f_1\otimes\dots\otimes
f_m)={\bf P}(f_1\otimes\dots\otimes f_m).\end{align*}

Thus, $\bf P$ is a bounded self-adjoint operator in ${\cal
A}({\cal H}_1,\dots,{\cal H}_l)$ satisfying ${\bf P}^2={\bf P}$,
and so $\bf P$ is an orthogonal projection. Hence, it remains only
to show that
\begin{equation}\label{763265}\Theta=\operatorname{Ker}{\bf
P}.\end{equation}

The inclusion $\Theta\subset\operatorname{Ker}{\bf P}$ follows
from \eqref{uuzfrdt}. Moreover, we have $${\cal A}({\cal
H}_1,\dots,{\cal H}_l)=\operatorname{Ker}{\bf
P}\oplus\operatorname{Im}{\bf P}.$$ Hence, to prove \eqref{763265}
it suffices to show that \begin{equation}\label{drdf}{\cal
A}({\cal H}_1,\dots,{\cal H}_l)=\Theta\oplus \operatorname{Im}{\bf
P}.\end{equation}

Let us fix arbitrary vectors $f_k\in{\cal H}_{i_k}$,
$i_k\in\{1,\dots,l\}$, $k=1,\dots,m$, $m\ge2$. We will now show
that the vector  $f_1\otimes\dots\otimes f_m$ can be represented
as a sum of vectors from $\Theta$ and $\operatorname{Im}{\bf P}$,
which will imply \eqref{drdf} (notice that ${\bf
P}\restriction{\cal A}_i({\cal H}_1,\dots,{\cal H}_l)={\bf 1}$,
$i=0,1$).

It is enough to show that, for each $\sigma\in S_m$, the vector
\begin{equation}\label{lawlio}F_\sigma:=(f_1\otimes\dots\otimes f_m)-
\sign(\sigma,i_1,\dots,i_m)( f_{\sigma(1)}\otimes\dots\otimes
f_{\sigma(m)})
\end{equation} belongs to $\Theta$, because \eqref{lawlio} yields $${\bf
P}(f_1\otimes\dots\otimes f_m)+\frac1{m!}\sum_{\sigma\in S_m
}F_\sigma= f_1\otimes\dots\otimes f_m.$$ But the inclusion
$F_\sigma\in\Theta$ can be proved by recurrent application of the
following  identity
\begin{align*}&\sign(\sigma,i_1,\dots,i_m)\left[(
f_{\sigma(1)}\otimes\dots\otimes f_{\sigma(m)})
-Q_{\sigma(s)}(f_{\sigma(1)}\otimes\dots\otimes
f_{\sigma(m)})\right]\\&\qquad
=\sign(\tau,i_1,\dots,i_m)(f_{\tau(1)}\otimes\dots\otimes
f_{\tau(m)}),\end{align*} where
 \begin{gather*}s:=\max\big\{\,r: r\in\{1,\dots,l\},\; \sigma(r)\ne r
 \,\big\},\\
 \tau(1,\dots,m):= (\sigma(1),\dots,\sigma(s-1),\sigma(s+1),\sigma(s),\sigma(s+2),\dots,
 \sigma(m)),\end{gather*} and by definition
\begin{gather*}Q_r (g_1\otimes\dots\otimes g_m):=
(g_1\otimes\dots\otimes g_m)-\\\text{}-
(-1)^{p(j_r)p(j_{r+1})}(g_1\otimes\dots\otimes g_{r-1}\otimes
g_{r+1}\otimes g_r\otimes g_{r+2}\otimes \dots\otimes g_m),\\
g_k\in{\cal H}_{j_k},\ j_k\in\{1,\dots,l\},\ k\in\{1,\dots,m\},\
r\in\{1,\dots,m-1\}.\end{gather*} Thus, Lemma~\ref{lemmas1} is
proven.

  Let us fix any orthonormal basis
$(e^{(i)}_k)_{k\ge1}$ in ${\cal H }_i$, $i=1,\dots,l$. Then, the
vectors $$ e^{(i_1)}_{k_1}\otimes\dots\otimes
e^{(i_m)}_{k_m},\qquad i_1,\dots,i_m\in\{1,\dots,l\},\
k_1,\dots,k_m\ge 1,$$ constitute an orthonormal basis in ${\cal
A}_m({\cal H}_1,\dots,{\cal H}_l)$, $m\in\N$. Therefore, by using
\eqref{ethess}, we conclude that  the following vectors constitute
an orthogonal basis  in ${\cal A}_{m,\, {\mathrm sym}}({\cal
H}_1,\dots,{\cal H}_l)$: $$ e^{(1)}_{k_1^{(1)}}\otimes\dots\otimes
e^{(1)}_{k^{(1)}_{r_1}}\otimes\dots\otimes
e^{(l)}_{k_1^{(l)}}\otimes\dots\otimes e^{(l)}_{k^{(l)}_{r_l}},
\qquad r_1,\dots,r_l\in\Z_+,\ r_1+\dots+r_l=m,$$ where
$k_1^{(i)}<k_2^{(i)}<\dots<k^{(i)}_{r_i}$ if $p(i)$ is odd, and
$k_1^{(i)}\le k_2^{(i)}\le \dots<k^{(i)}_{r_i}$ if $p(i)$ is even.
For any such  vector, we get \begin{gather}\big \|{\bf P}(
e^{(1)}_{k_1^{(1)}}\otimes\dots\otimes
e^{(1)}_{k^{(1)}_{r_1}}\otimes\dots\otimes
e^{(l)}_{k_1^{(l)}}\otimes\dots\otimes
e^{(l)}_{k^{(l)}_{r_l}})\big\|^2=\notag\\= \big({\bf P}(
e^{(1)}_{k_1^{(1)}}\otimes\dots\otimes
e^{(1)}_{k^{(1)}_{r_1}}\otimes\dots\otimes
e^{(l)}_{k_1^{(l)}}\otimes\dots\otimes
e^{(l)}_{k^{(l)}_{r_l}}),e^{(1)}_{k_1^{(1)}}\otimes\dots\otimes
e^{(1)}_{k^{(1)}_{r_1}}\otimes\dots\otimes
e^{(l)}_{k_1^{(l)}}\otimes\dots\otimes
e^{(l)}_{k^{(l)}_{r_l}}\big)\notag\\=\frac1{m!}\,\prod_{j=1}^l\bigg[
\sum_{\sigma_j\in
S_{r_j}}\sign(\sigma_j,\underbrace{j,\dots,j}_{\text{$r_j$
times}})\, \big(e^{(j)}_{k^{(j)}_{\sigma_j(1)}}\otimes\dots\otimes
e^{(j)}_{k^{(j)}_{\sigma_j(r_j)}},
e^{(j)}_{k_1^{(j)}}\otimes\dots\otimes
e^{(j)}_{k^{(j)}_{r_j}}\big)\bigg]\notag\\=\frac1{m!}\,\prod_{j=1}^l\bigg[
\sum_{\sigma_j\in S_{r_j}}{\frak
S}(\sigma_j,j)\big(e^{(j)}_{k^{(j)}_{\sigma_j(1)}}\otimes\dots\otimes
e^{(j)}_{k^{(j)}_{\sigma_j(r_j)}},
e^{(j)}_{k_1^{(j)}}\otimes\dots\otimes
e^{(j)}_{k^{(j)}_{r_j}}\big)\bigg]\notag\\ =\frac{r_1!\dotsm
r_l!}{m!}\, \prod _{j=1}^l\big\|
e^{(j)}_{k^{(j)}_1}\overset{p(j)}{\diamond}\dotsm
\overset{p(j)}{\diamond}e^{(j)}_{k^{(j)}_{r_j}}\|^2_{ {\cal
H}_j^{\overset{p(j)}{\diamond} r_j }     },\label{klawli}
\end{gather}
where \begin{equation} {\frak S}(\sigma_j,j)=\begin{cases}\sign
\sigma_j,&\text{if $p(j)$ is odd,}\\ 1,&\text{if $p(j)$ is
even.}\end{cases}\label{johjft}\end{equation} From \eqref{klawli}
and \eqref{johjft} the conclusion of Lemma~\ref {lemmas2}
trivially follows.

\end{document}